\documentclass[12pt]{amsart}
\usepackage[all]{xy}
\usepackage{amssymb}

\title{On the class of Benson's cofibrant modules}
\author{Ioannis Emmanouil and Wei Ren}

\newtheorem{Lemma}{Lemma}[section]
\newtheorem{Proposition}[Lemma]{Proposition}
\newtheorem{Theorem}[Lemma]{Theorem}
\newtheorem{Corollary}[Lemma]{Corollary}

\begin{document}

\begin{abstract}
In this paper, we study the class of cofibrant modules over a
group algebra $kG$, that were introduced by Benson. We prove
that this class is always the left-hand side of a complete
hereditary and projective cotorsion pair. We also examine the
relation between cofibrant and Gorenstein projective modules
and the behaviour of cofibrant modules with respect to the
induction functor from subgroups. As an application, we show
that the cofibrant cotorsion pair induces a monoidal model
structure on the category of $kG$-modules over {\em any}
commutative coefficient ring $k$, provided that the group
$G$ is obtained from the class of finite groups by iterated
applications of Kropholler's operation ${\scriptstyle{{\bf LH}}}$
and Talelli's operation $\Phi$. The corresponding symmetric
monoidal homotopy category is equivalent to the stable category
of cofibrant $kG$-modules.
\end{abstract}

\thanks{}

\maketitle
\tableofcontents

\addtocounter{section}{-1}
\section{Introduction}

\noindent
Gorenstein homological algebra is the relative homological theory,
which is based upon the classes of Gorenstein projective, Gorenstein
injective and Gorenstein flat modules. The theory has experienced a
rapid development during the past several years and has found interesting
applications in the representation theory of Artin algebras. Gorenstein
homological algebra has found applications in cohomological group
theory as well. For example, if $G$ is a group and $\underline{E}G$
is a model for the classifying space for proper actions of $G$, then
the augmented cellular chain complex of $\underline{E}G$ is a resolution
of the trivial $\mathbb{Z}G$-module $\mathbb{Z}$ by Gorenstein projective
$\mathbb{Z}G$-modules. In this way, the Gorenstein cohomological dimension
of groups over $\mathbb{Z}$ is proposed in \cite{BDT} to serve as an
algebraic invariant whose finiteness characterizes those groups that
admit a finite dimensional model for the classifying space for proper
actions.

Building upon the relation between the class of finite groups and
the class of those groups that admit a cellular action on a finite
dimensional contractible CW-complex with finite cell stabilizers,
Kropholler defined in \cite{Kro} the class
${\scriptstyle{{\bf LH}}}\mathfrak{F}$, that consists of those groups
whose finitely generated subgroups admit a hierarchical description
starting form the class of finite groups. In fact, Kropholler's
operation ${\scriptstyle{{\bf LH}}}$ may be performed to any class
of groups. Another example of such an operation on group classes is
due to Talelli, who introduced the groups of type $\Phi$ in \cite{T},
in order to construct examples that admit a finite dimensional model
for the classifying space for proper actions. The definition of
these groups is based on the attempt to describe modules of finite
projective dimension over the group, in terms of the finiteness of
the projective dimension of their restriction to finite subgroups.
Talelli's operation $\Phi$ may be performed with base any class
of groups as well. We consider the class $\overline{\mathfrak{F}}$,
that is obtained from the class $\mathfrak{F}$ of finite groups by
a transfinite iteration of the operations ${\scriptstyle{{\bf LH}}}$
and $\Phi$; the precise definition of $\overline{\mathfrak{F}}$ is
given in $\S $1.II below. This class is very big and it does require
some effort to construct groups which are not contained therein; see
the appendix of \cite{ET4} though. Often, homological properties of
finite groups can be extended to $\overline{\mathfrak{F}}$-groups;
examples of this phenomenon can be found in \cite{ET3}, \cite{ER},
\cite{ET4}.

A useful tool for the construction of examples of Gorenstein
projective $\mathbb{Z}G$-modules is the module $B(G,\mathbb{Z})$
of bounded functions on $G$ with values in $\mathbb{Z}$, which
is introduced in \cite{KT}. Cornick and Kropholler proved in
\cite{CK1} that the class of Gorenstein projective
$\mathbb{Z}G$-modules contains any $\mathbb{Z}G$-module $M$,
for which the diagonal $\mathbb{Z}G$-module
$M \otimes_{\mathbb{Z}} B(G,\mathbb{Z})$ is projective. Benson
\cite{Ben} realized the importance of the latter condition for
the construction of model structures in module categories over
the group algebras of certain groups and termed the
$\mathbb{Z}G$-modules satisfying that condition as cofibrant.
It is shown by Dembegioti and Talelli in \cite{DT} that any
Gorenstein projective $\mathbb{Z}G$-module is necessarily
cofibrant if $G$ is an ${\scriptstyle{{\bf LH}}}\mathfrak{F}$-group.
The equality between the classes ${\tt Cof}(\mathbb{Z}G)$ and
${\tt GProj}(\mathbb{Z}G)$ of cofibrant and Gorenstein projective
$\mathbb{Z}G$-modules, which is conjectured to hold over any group
$G$ in [loc.cit.], is an interesting problem in cohomological group
theory. If the equality between these classes is indeed true, then
the study of Gorenstein projective $\mathbb{Z}G$-modules would
become considerably simpler; it would be essentially reduced to
certain questions involving classical homological algebra notions.
For example, it would follow that being Gorenstein projective is
a module property, that is closed under restriction to subgroups.
More generally, the question as to whether the classes
${\tt Cof}(kG)$ and ${\tt GProj}(kG)$ of cofibrant and Gorenstein
projective $kG$-modules are equal has been studied for any
commutative coefficient ring $k$. As an indicative result in this
direction, we note that the two module classes are equal if $G$ is
an ${\scriptstyle{{\bf LH}}}\mathfrak{F}$-group and $k$ has finite
global dimension; see \cite{Bi}. In this paper, we show the following
result (cf.\ Corollary 4.8):

\medskip

\noindent
{\bf Theorem A.}
{\em Assume that any acyclic complex of projective $k$-modules
is contractible. If $G$ is an $\overline{\mathfrak{F}}$-group,
then ${\tt Cof}(kG) = {\tt GProj}(kG)$.}

\medskip

\noindent
We note that the hypothesis made on $k$ is satisfied if $k$ has
finite weak global dimension. On the other hand, there are examples
of commutative rings $k$ (that do not satisfy the hypothesis in
Theorem A), over which ${\tt Cof}(kG) \neq {\tt GProj}(kG)$ for
any group $G$; cf.\ Remark 4.1(iii).

It is shown by Cort\'{e}s-Izurdiaga and \v{S}aroch in \cite{CIS}
that the pair $\left( {\tt GProj}(R),{\tt GProj}(R)^{\perp} \right)$
is a hereditary cotorsion pair in the category of $R$-modules over
any ring $R$. An interesting problem in Gorenstein homological
algebra is to examine whether this cotorsion pair is complete. The
importance of completeness for a cotorsion pair (i.e.\ the importance
of the existence of suitable approximations) in the development of
the associated relative homological theory is analysed by Auslander
and Reiten in \cite{AR}. Our main result regarding cofibrant modules
in this paper is the following (cf.\ Theorem 2.4 and Proposition 3.3):

\medskip

\noindent
{\bf Theorem B.}
{\em Assume that $k$ is any commutative ring and let $G$ be any
group. Then, the pair
$\left( {\tt Cof}(kG),{\tt Cof}(kG)^{\perp} \right)$ is a complete
hereditary cotorsion pair in the category of $kG$-modules.}

\medskip

\noindent
The stable category of modules over the algebra of a finite group
$G$ with coefficients in a field $F$ is a well-known example of
a compactly generated tensor triangulated category. This stable
category is a central theme of study in (modular) representation
theory. We note that any $FG$-module is cofibrant in this case.
In general, if $k$ is any commutative coefficient ring and $G$
 is any group, the exact subcategory of cofibrant $kG$-modules
is Frobenius with projective-injective objects the projective
$kG$-modules and hence we may define the associated stable
category $\underline{\tt Cof}(kG)$. We shall prove the following
result in Section 7; cf.\ Corollary 7.5.

\medskip

\noindent
{\bf Theorem C.}
{\em Assume that $k$ is any commutative ring and let $G$ be
an $\overline{\mathfrak{F}}$-group. Then, the tensor product
of $kG$-modules (with the diagonal action of $G$) endows the
stable category $\underline{\tt Cof}(kG)$ with the structure
of a symmetric monoidal category.}

\medskip

\noindent
The existence of a tensor unit is not immediately clear, since
the trivial $kG$-module $k$ is only rarely cofibrant.\footnote{If
$G$ has finite cohomological dimension over $k$, then the trivial
$kG$-module $k$ is cofibrant only if $G$ is finite.} It will turn
out that a tensor unit is provided by a special
${\tt Cof}(kG)$-precover of $k$. In fact, the stable category of
cofibrant modules is equivalent as a triangulated category to the
homotopy category of the model structure on the category of all
$kG$-modules, that is induced by the cotorsion pair of Theorem B.
The key fact is that this model structure is monoidal with respect
to the tensor product of $kG$-modules, in the sense of
\cite[Definition 4.2.6]{Hov1}, for any $\overline{\mathfrak{F}}$-group
$G$.

\medskip

\noindent
{\em Notations and terminology.}
Unless otherwise specified, we work with modules over the
group algebra $kG$ of a group $G$ with coefficients in a
unital commutative ring $k$.

\section{Preliminaries}

\noindent
In this section, we record certain prerequisite notions that are
used throughout the paper. These notions concern basic definitions
regarding modules over group algebras, hierarchically defined group
classes, cotorsion pairs and model categories.

\vspace{0.1in}

\noindent
{\sc I.\ Modules over group algebras.}
We consider a unital commutative ring $k$ and let $G$ be a group.
Using the diagonal action of $G$, the tensor product
$M \otimes_k N$ of two $kG$-modules $M,N$ is also a $kG$-module;
we let $g \cdot (x \otimes y) = gx \otimes gy \in M \otimes_k N$
for any $g \in G$, $x \in M$ and $y \in N$. If the $kG$-module
$M$ is projective and $N$ is projective as a $k$-module, then the
diagonal $kG$-module $M \otimes_k N$ is projective as well; cf.\
\cite[Chapter III, Corollary 5.7]{Bro}. We also note that for any two
$kG$-modules $M,N$ the $k$-module ${\rm Hom}_k(N,M)$ admits the
structure of a $kG$-module, with the group $G$ acting diagonally;
we now let $(g \cdot f)(x) = gf(g^{-1}x)\in M$ for any $g \in G$,
$f \in {\rm Hom}_k(N,M)$ and $x \in N$. For any three $kG$-modules
$L,M,N$ the natural isomorphism of $k$-modules
${\rm Hom}_k(L \otimes_k N,M) \cong {\rm Hom}_k(L, {\rm Hom}_k(N,M))$
is actually an isomorphism of $kG$-modules, where $G$ acts diagonally
on the tensor product and the Hom-groups.

Let $B(G,\mathbb{Z})$ be the $\mathbb{Z}G$-module consisting
of all bounded functions from $G$ to $\mathbb{Z}$, introduced
in \cite{KT}. The $\mathbb{Z}G$-module $B(G,\mathbb{Z})$ is
$\mathbb{Z}$-free; in fact, $B(G,\mathbb{Z})$ is $\mathbb{Z}H$-free
for any finite subgroup $H \subseteq G$. If $n$ is an integer, then
the constant function $\iota(n) \in B(G,\mathbb{Z})$ with value $n$
is invariant under the action of $G$. The map
$\iota : \mathbb{Z} \longrightarrow B(G,\mathbb{Z})$ which is
defined in this way is therefore $\mathbb{Z}G$-linear. Moreover,
$\iota$ is a $\mathbb{Z}$-split monomorphism; an additive
splitting for $\iota$ may be obtained by evaluating functions
at the identity element of $G$. The $kG$-module
$B(G,k) = B(G,\mathbb{Z}) \otimes_{\mathbb{Z}}k$ may be
identified with the $kG$-module of all functions from $G$
to $k$ that admit finitely many values. It is free as a
$kH$-module for any finite subgroup $H \subseteq G$. We also
note that $\iota$ induces a $k$-split $kG$-linear monomorphism
$\iota \otimes 1 : k \longrightarrow B(G,k)$. Finally, if
$G_0 \subseteq G$ is a subgroup, then $B(G_0,k)$ is a direct
summand of the restricted $kG_0$-module $\mbox{res}_{G_0}^GB(G,k)$.
For simplicity of notation, if the coefficient ring $k$ and the group
$G$ are understood, we shall denote $B(G,k)$ by $B(G)$ or $B$ and let
$\overline{B} = \mbox{coker} \, (\iota \otimes 1)$.

\vspace{0.1in}

\noindent
{\sc II.\ Operations and the closure of group classes.}
An operation $\mathcal{O}$ on group classes assigns a group class
$\mathcal{O}\mathfrak{C}$ to each group class $\mathfrak{C}$, in
such a way that:
\newline
(i) $\mathfrak{C} \subseteq \mathcal{O}\mathfrak{C}$ for any group
class $\mathfrak{C}$ and
\newline
(ii) $\mathcal{O}\mathfrak{C}_1 \subseteq \mathcal{O}\mathfrak{C}_2$
for any two group classes $\mathfrak{C}_1,\mathfrak{C}_2$ with
$\mathfrak{C}_1 \subseteq \mathfrak{C}_2$.
\newline
The operation is continuous if it also has the following property:
\newline
(iii) Assume that $\mathfrak{C}_{\alpha}$ is a group class defined
for each ordinal number $\alpha$, such that
$\mathfrak{C}_{\alpha} \subseteq \mathfrak{C}_{\beta}$ for all
$\alpha \leq \beta$, and let $\mathfrak{C}$ be the class consisting
of those groups which are contained in $\mathfrak{C}_{\alpha}$ for
some ordinal number $\alpha$. Then, the group class
$\mathcal{O}\mathfrak{C}$ is the class consisting of those groups
which are contained in $\mathcal{O}\mathfrak{C}_{\alpha}$ for some
$\alpha$.

A group class $\mathfrak{C}$ is said to be closed under an operation
$\mathcal{O}$ if $\mathcal{O}\mathfrak{C} = \mathfrak{C}$. The
$\mathcal{O}$-closure $\overline{\mathfrak{C}}$ of a group class
$\mathfrak{C}$ is the smallest $\mathcal{O}$-closed group class
containing $\mathfrak{C}$. If the operation $\mathcal{O}$ is
continuous, we may describe the $\overline{\mathfrak{C}}$-groups
by a hierarchical process, as follows:  We define the classes
$\mathfrak{C}_{\alpha}$ for any ordinal number $\alpha$, by
letting $\mathfrak{C}_0 = \mathfrak{C}$,
$\mathfrak{C}_{\alpha +1} = \mathcal{O}\mathfrak{C}_{\alpha}$
for any ordinal $\alpha$ and
$\mathfrak{C}_{\alpha} =
 \bigcup_{\beta < \alpha}\mathfrak{C}_{\beta}$
for any limit ordinal $\alpha$. Then, $\overline{\mathfrak{C}}$
is the class consisting of those groups which are contained in
$\mathfrak{C}_{\alpha}$, for some ordinal $\alpha$.

If $\mathcal{O}_1, \mathcal{O}_2, \ldots , \mathcal{O}_n$ are
(continuous) operations on group classes, then we may consider
the  (continuous) operation
$\mathcal{O}  = \mathcal{O}_1 \cup \mathcal{O}_2 \cup \ldots
 \cup \mathcal{O}_n$,
which is defined by letting
$\mathcal{O}\mathfrak{C} = \mathcal{O}_1\mathfrak{C} \cup
 \mathcal{O}_2\mathfrak{C} \cup \ldots \cup
 \mathcal{O}_n\mathfrak{C}$
for any group class $\mathfrak{C}$. The $\mathcal{O}$-closure
of $\mathfrak{C}$ is the smallest group class which contains
$\mathfrak{C}$ and is closed under all of the operations
$\mathcal{O}_1, \mathcal{O}_2, \ldots , \mathcal{O}_n$.

We are interested in the continuous operation
${\scriptstyle{{\bf LH}}}$ that was introduced by Kropholler
in \cite{Kro}. If $\mathfrak{C}$ is any group class, then we
define for each ordinal $\alpha$ the group class
${\scriptstyle{{\bf H}}}_{\alpha}\mathfrak{C}$, using transfinite
induction, as follows: We let
${\scriptstyle{{\bf H}}}_0\mathfrak{C} = \mathfrak{C}$. For an
ordinal $\alpha >0$, the class
${\scriptstyle{{\bf H}}}_{\alpha}\mathfrak{C}$ consists of those
groups $G$ which admit a cellular action on a finite dimensional
contractible CW-complex $X$, in such a way that each isotropy
subgroup of the action belongs to
${\scriptstyle{{\bf H}}}_{\beta}\mathfrak{C}$ for some ordinal
$\beta < \alpha$ (that may depend upon the particular cell). If
$G,X$ are as above, $k$ is a commutative ring and $M$ is a $kG$-module,
the cellular chain complex of $X$ induces an exact sequences
of $kG$-modules
\[ 0 \longrightarrow M_d \longrightarrow \cdots
     \longrightarrow M_1 \longrightarrow M_0
     \longrightarrow M \longrightarrow 0 , \]
where $d$ is the geometric dimension of $X$ and the $M_i$'s are
direct sums of $kG$-modules of the form $\mbox{ind}_H^G\mbox{res}_H^GM$,
for some ${\scriptstyle{{\bf H}}}_{\beta}\mathfrak{C}$-subgroup $H$
of $G$ with $\beta < \alpha$. We say that a group belongs to
${\scriptstyle{{\bf H}}}\mathfrak{C}$ if it belongs to
${\scriptstyle{{\bf H}}}_{\alpha}\mathfrak{C}$ for some $\alpha$.
Then, the class ${\scriptstyle{{\bf LH}}}\mathfrak{C}$ consists
of those groups $G$, all of whose finitely generated subgroups
$H$ are contained in some
${\scriptstyle{{\bf H}}}\mathfrak{C}$-subgroup $K=K(H) \subseteq G$.
The continuity of the operation ${\scriptstyle{{\bf LH}}}$ is
proved in \cite[Lemma 1.3]{ET3}. If the class $\mathfrak{C}$
is subgroup-closed, then the classes
${\scriptstyle{{\bf H}}}\mathfrak{C}$ and
${\scriptstyle{{\bf LH}}}\mathfrak{C}$ are subgroup-closed as
well; in this case, ${\scriptstyle{{\bf LH}}}\mathfrak{C}$
consists of those groups all of whose finitely generated
subgroups are ${\scriptstyle{{\bf H}}}\mathfrak{C}$-groups.

The concept of groups of type $\Phi$, that were introduced by
Talelli in \cite{T}, corresponds to yet another continuous
operation, applied to the class $\mathfrak{F}$ of finite groups.
More generally, if $k$ is a commutative ring and $\mathfrak{C}$
is a class of groups, then we define $\Phi\mathfrak{C}$ as the
class consisting of those groups $G$, for which the following
two conditions are equivalent for any $kG$-module $M$:
\newline
(i) $\mbox{pd}_{kG}M < \infty$,
\newline
(ii) there is an integer $n$, such that
$\mbox{pd}_{kH}\mbox{res}_H^GM \leq n$ for any
$\mathfrak{C}$-subgroup $H \subseteq G$.
\newline
See also \cite{Bi} and \cite{MS}. The continuity of $\Phi$
is proved in \cite[Lemma 1.4]{ET3}. We also consider the
version of the operation $\Phi$ for the flat dimension of
modules and define for any group class $\mathfrak{C}$ the
class $\Phi_{flat} \mathfrak{C}$ accordingly.

\vspace{0.1in}

\noindent
{\sc III.\ Cotorsion pairs.}
The ${\rm Ext}^1$-pairing induces an orthogonality relation
between modules over a ring $R$. If {\tt S} is a class of
modules, then the left orthogonal $^{\perp}{\tt S}$ of {\tt S}
is the class consisting of those modules $M$, for which
${\rm Ext}^1_R(M,S)=0$ for all $S \in {\tt S}$. The right
orthogonal ${\tt S}^{\perp}$ of {\tt S} is the class consisting
of those modules $N$, for which ${\rm Ext}^1_R(S,N)=0$ for all
$S \in {\tt S}$. Any left orthogonal class is closed under
transfinite extensions, in the sense that we shall now explain.

Let {\tt T} be any class of modules. A continuous ascending
filtration of a module $M$ by modules in {\tt T} is a family
of submodules $(M_{\alpha})_{\alpha < \gamma}$ of $M$, which
is indexed by an ordinal $\gamma$, such that:

(i) $M_0=0$ and $M = \bigcup_{\alpha < \gamma} M_{\alpha}$,

(ii) $M_{\alpha} \subseteq M_{\beta}$ whenever $\alpha$ and
$\beta$ are two ordinals with $\alpha < \beta < \gamma$,

(iii) $M_{\beta} = \bigcup_{\alpha < \beta} M_{\alpha}$ for any
limit ordinal $\beta < \gamma$ and

(iv) $M_{\alpha +1}/M_{\alpha} \in {\tt T}$ for any ordinal
$\alpha$ with $\alpha + 1 < \gamma$.
\newline
We denote by Filt-{\tt T} the class of those modules which admit
a continuous ascending filtration by modules in {\tt T}. Then,
Eklof's lemma \cite[Theorem 7.3.4]{EJ} asserts that
Filt-$\left( ^{\perp}{\tt S} \right) \! \subseteq \! ^{\perp}{\tt S}$
for any class of modules {\tt S}.

If ${\tt U}$ and ${\tt V}$ are two classes of modules, then
we say that the pair $({\tt U},{\tt V})$ is a cotorsion pair
in the category of modules (cf.\ \cite[Definition 7.1.2]{EJ})
if ${\tt U} = {^{\perp}{\tt V}}$ and ${\tt U} ^{\perp} = {\tt V}$.
In that case, the class {\tt U} is closed under extensions,
direct summands and (in view of Eklof's lemma) continuous
ascending filtrations. The cotorsion pair $({\tt U},{\tt V})$
is called hereditary if ${\rm Ext}^i_R(U,V)=0$ for all $i>0$
and all modules $U \in {\tt U}$ and $V \in {\tt V}$. Equivalently,
the cotorsion pair is hereditary if {\tt U} is closed under
kernels of epimorphisms or else if {\tt V} is closed under
cokernels of monomorphisms. The cotorsion pair is complete
if for any module $M$ there exist short exact sequences
\[ 0 \longrightarrow V \longrightarrow U \longrightarrow M
     \longrightarrow 0
   \;\;\; \mbox{and } \;\;\;
   0 \longrightarrow M \longrightarrow V' \longrightarrow U'
     \longrightarrow 0 , \]
where $U,U' \in {\tt U}$ and $V,V' \in {\tt V}$. The existence
of these short exact sequences is also referred to by saying that
any module admits a special {\tt U}-precover and a special
{\tt V}-preenvelope respectively. The kernel of the cotorsion
pair $({\tt U},{\tt V})$ is the intersection
${\tt U} \cap {\tt V}$; we say that the cotorsion pair is
projective if its kernel coincides with the class of projective
modules. The cotorsion pair $({\tt U},{\tt V})$ is cogenerated
by a set of modules {\tt S} if ${\tt V} = { \tt S} ^{\perp}$, so
that ${\tt U} = \, \! ^{\perp} \! \left( {\tt S} ^{\perp} \right)$.
In that case, the left-hand class {\tt U} of the pair consists
precisely of the direct summands of Filt-$({\tt S} \cup \{ R \})$
modules; cf.\ \cite[Corollary 7.3.5]{EJ}. Eklof and Trlifaj have
established in \cite[Theorem 10]{ETz} the completeness of any
cotorsion pair which is cogenerated by a set of modules.

\vspace{0.1in}

\noindent
{\sc IV.\ Hovey triples and model categories.}
A Hovey triple in the category of modules over a ring $R$ is a
triple $({\mathcal C} , {\mathcal W} , {\mathcal F})$ of module
classes, which are such that the pairs
$({\mathcal C} , {\mathcal W} \cap {\mathcal F})$ and
$({\mathcal C} \cap {\mathcal W} , {\mathcal F})$ are complete
cotorsion pairs and the class ${\mathcal W}$ is thick (i.e.\
$\mathcal{W}$ is closed under direct summands and satisfies
the 2-out-of-3 property for short exact sequences). The work
of Hovey \cite{Hov2} and Gillespie \cite{Gil} provides a
bijection between Hovey triples and certain model structures
in the module category $R$-Mod; cf.\ \cite[Theorem 3.3]{Gil}.
A model structure $\mathcal{M}$ on a category consists in
specifying three classes of morphisms, the fibrations, the
cofibrations and the weak equivalences, that have certain
properties which are reminiscent of the usual properties of
the corresponding classes of maps in topology. The homotopy
category $\mbox{Ho}(\mathcal{M})$ of $\mathcal{M}$, which is
obtained by formally inverting the weak equivalences, is
triangulated. A comprehensive account of model categories,
a notion originally introduced by Quillen in \cite{Q}, may
be found in Hovey's book \cite{Hov1}. As shown in
\cite[Proposition 5.2]{Gil}, for a model structure
$\mathcal{M}$ in $R$-Mod that corresponds to a Hovey triple
$({\mathcal C} , {\mathcal W} , {\mathcal F})$ that induces
two {\em hereditary} complete cotorsion pairs, the class
${\mathcal C} \cap {\mathcal F}$ is a Frobenius category
with projective-injective objects the modules in
${\mathcal C} \cap {\mathcal W} \cap {\mathcal F}$. As shown
by Happel \cite{Ha}, the stable category
$\underline{{\mathcal C} \cap {\mathcal F}}$, i.e.\
${\mathcal C} \cap {\mathcal F}$ modulo its projective-injective
objects ${\mathcal C} \cap {\mathcal W} \cap {\mathcal F}$, is
triangulated. Using \cite[Proposition 4.4 and Corollary 4.8]{Gil},
it follows that the homotopy category $\mbox{Ho}(\mathcal{M})$
of such a model structure $\mathcal{M}$ is triangulated equivalent
to $\underline{{\mathcal C} \cap {\mathcal F}}$.

\section{Completeness of the cotorsion pair $\left( {\tt Cof}(kG),{\tt Cof}(kG)^{\perp} \right)$}

\noindent
Following Benson \cite{Ben}, we say that a $kG$-module $M$ is cofibrant
if the (diagonal) $kG$-module $M \otimes_k B$ is projective; we denote
by ${\tt Cof}(kG)$ the class of cofibrant modules. In this section, we
show that the class of cofibrant modules is a Kaplansky class and conclude
that $({\tt Cof}(kG),{\tt Cof}(kG)^{\perp})$ is a cotorsion pair, which
is cogenerated by a set of modules.

We begin with an auxiliary simple lemma.

\begin{Lemma}
The class ${\tt Cof}(kG)$ is closed under direct summands and
transfinite extensions.
\end{Lemma}

\begin{proof}
If the $kG$-module $M$ is a direct summand of a cofibrant
module $M'$, then $M \otimes_k B$ is a direct summand of the
projective $kG$-module $M' \otimes_k B$. It follows that the
$kG$-module $M \otimes_k B$ is itself projective, so that $M$
is cofibrant. Hence, ${\tt Cof}(kG)$ is closed under direct
summands.

Let $(N_{\alpha})_{\alpha < \tau}$ be a continuous
ascending filtration of a $kG$-module $N$ by cofibrant modules.
Then, $(N_{\alpha} \otimes_k B)_{\alpha < \tau}$ is a continuous
ascending filtration of the $kG$-module $N \otimes _k B$ by
projective $kG$-modules. Since the class ${\tt Proj}(kG)$ of
projective $kG$-modules is closed under transfinite extensions,
it follows that $N \otimes_k B \in {\tt Proj}(kG)$, so that $N$
is cofibrant.
\end{proof}

\noindent
We recall from \cite[Definition 2.1]{ELR} that a class $\mathcal{K}$
of $kG$-modules is called a Kaplansky class if there exists a cardinal
number $\lambda$, such that for any module $M \in \mathcal{K}$ and any
element $x \in M$ there exists a submodule $N \subseteq M$ with $x \in N$,
$\mbox{card} \, N \leq \lambda$ and $N,M/N \in \mathcal{K}$. Given a
cardinal number $\mu$, we say that a $kG$-module is $\mu$-generated if
it can be generated by a subset of cardinality $\leq \mu$.

\begin{Proposition}\label{prop:KapC}
The class ${\tt Cof}(kG)$ is a Kaplansky class.
\end{Proposition}

\begin{proof}
Let $\kappa = \mbox{card} \, k$, $\gamma = \mbox{card} \, G$
and $\lambda = \max \{ \aleph_0, \kappa, 2^{\gamma} \}$. Then,
$\lambda$ is an infinite cardinal, which is greater than or equal
to both $\mbox{card} \ kG$ and $\mbox{card} \, B$.\footnote{The
inequality $\mbox{card} \, B \leq \lambda$ follows since any
element of $B$ is a $k$-linear combination of suitable characteristic
functions on subsets of $G$.} Moreover, it is easily seen that for
any $\lambda$-generated $kG$-module $K$ (i.e.\ for any $kG$-module
$K$ with $\mbox{card} \, K \leq \lambda$) the $kG$-module
$K \otimes_k B$ is also $\lambda$-generated (i.e.\
$\mbox{card} \, (K \otimes_k B) \leq \lambda$). We prove that
this choice of $\lambda$ works to establish that ${\tt Cof}(kG)$
is a Kaplansky class.

To that end, we consider a cofibrant module $M$ and fix an element
$x \in M$. Since $M \otimes_k B$ is a projective $kG$-module, a
result by Kaplansky \cite{Kap} implies that it can be expressed as the
direct sum of a family of countably generated projective
$kG$-modules. We fix such a decomposition and write
$M \otimes_k B = \bigoplus_{i \in I}P_i$, for a suitable family
$(P_i)_{i \in I}$ of countably generated projective $kG$-modules.
For any subset $J \subseteq I$ we define
$P(J) = \bigoplus_{i \in J}P_i$; in particular, $P(\emptyset)=0$
and $P(I) = M \otimes_k B$. Using induction, we construct a
sequence of pairs $(I_n,M_n)_n$, such that:

(i) $(I_n)_n$ is an increasing sequence of subsets of $I$, with
$\mbox{card} \, I_n \leq \lambda$ for all $n \geq 0$,

(ii) $(M_n)_n$ is an increasing sequence of $\lambda$-generated
$kG$-submodules of $M$,

(iii) $I_0=\emptyset$ and $M_0 = kG \cdot x$,

(iv) $P(I_n) \subseteq M_n \otimes_k B \subseteq P(I_{n+1})
             \subseteq P(I) = M \otimes_k B$
for all $n \geq 0$.
\newline
Of course, we begin the induction by letting $I_0=\emptyset$ and
$M_0 = kG \cdot x$. Having defined the pairs
$(I_0,M_0),(I_1,M_1), \ldots , (I_n,M_n)$, we proceed with the
inductive step and define $(I_{n+1},M_{n+1})$ as follows: Since
the $kG$-module $M_n$ is $\lambda$-generated, the $kG$-module
$M_n \otimes_k B$ is also $\lambda$-generated. Expressing each
generator of
$M_n \otimes_k B \subseteq M \otimes_k B = \bigoplus_{i \in I}P_i$
as a finite sum of elements of the $P_i$'s, we conclude that there
exists a subset $I_{n+1} \subseteq I$ which contains $I_n$ and has
cardinality $\leq \lambda$, such that
$M_n \otimes_k B \subseteq P(I_{n+1})$. Since
$\mbox{card} \, I_{n+1} \leq \lambda$ and the $kG$-module $P_i$
is countably generated for all $i \in I$, the $kG$-module
$P(I_{n+1})$ is $\lambda$-generated. Expressing each generator of
$P(I_{n+1}) \subseteq P(I) = M \otimes_k B$ as a finite sum of
elementary tensors, we conclude that there is a $\lambda$-generated
$kG$-submodule $M_{n+1} \subseteq M$ containing $M_n$, such
that $P(I_{n+1}) \subseteq M_{n+1} \otimes_k B$. This completes
the inductive step of the construction.

Having constructed the sequence $(I_n,M_n)_n$ with properties
(i)-(iv) above, we let $I' = \bigcup_nI_n$ and $M' = \bigcup_nM_n$.
Then, $I' \subseteq I$ has cardinality $\leq \lambda$ and $M'$
is a $\lambda$-generated $kG$-submodule of $M$, which contains
$M_0$ (and hence $x$), such that
\[ M' \otimes_k B =
   {\textstyle{\bigcup}_n} (M_n \otimes_k B) =
   {\textstyle{\bigcup}_n} P(I_n) =
   P(I') \subseteq P(I) = M \otimes_k B . \]
Here, the second equality follows from property (iv). In
particular, the $kG$-module $M' \otimes_k B$ is projective
and hence $M'$ is cofibrant. The quotient module $M/M'$ is
also cofibrant, since
\[ (M/M') \otimes_k B = (M \otimes_k B)/(M' \otimes_k B) =
   P(I)/P(I') = P(I \setminus I') \]
is a projective $kG$-module. We have therefore completed the
proof that ${\tt Cof}(kG)$ is a Kaplansky class.
\end{proof}

\noindent
We shall use Proposition \ref{prop:KapC} (and its proof) to construct
filtrations of cofibrant modules with successive quotients
cofibrant modules of controlled size.

\begin{Proposition}\label{prop:caf}
There is a cardinal number $\lambda$, such that any cofibrant
module admits a continuous ascending filtration by
$\lambda$-generated cofibrant modules.
\end{Proposition}

\begin{proof}
Let $\lambda$ be the cardinal number defined in the proof
of Proposition \ref{prop:KapC}. We consider a cofibrant module $M$ and
fix a decomposition $M \otimes_k B = \bigoplus_{i \in I}P_i$
of the projective $kG$-module $M \otimes_k B$ into the direct
sum of a family $(P_i)_{i \in I}$ of countably generated
projective $kG$-modules, as in the proof of Proposition \ref{prop:KapC}.
As before, for any subset $J \subseteq I$ we define
$P(J) = \bigoplus_{i \in J}P_i$. We also choose a well-ordered
set of generators $(x_{\alpha})_{\alpha < \tau}$ of the
$kG$-module $M$, where $\tau$ is a suitable cardinal number.
It only suffices to consider the case where $\tau > \lambda$,
since otherwise $M$ is already $\lambda$-generated. Using
transfinite induction, we construct a family of pairs
$(I_{\alpha},M_{\alpha})_{\alpha \leq \tau}$, such that:

(i) $(I_{\alpha})_{\alpha \leq \tau}$ is a continuous ascending
family of subsets of $I$,

(ii) $\mbox{card} \, (I_{\alpha +1} \setminus I_{\alpha})
      \leq \lambda$
for any $\alpha < \tau$,

(iii) $(M_{\alpha})_{\alpha \leq \tau}$ is a continuous ascending
family of $kG$-submodules of $M$,

(iv) the $kG$-module $M_{\alpha +1}/M_{\alpha}$ is
$\lambda$-generated for any $\alpha < \tau$,

(v) $x_{\alpha} \in M_{\alpha +1}$ for any $\alpha < \tau$,

(vi) $I_0 = \emptyset$ and $M_0=0$,

(vii) $M_{\alpha} \otimes_k B = P(I_{\alpha}) \subseteq
       P(I) = M \otimes_k B$
for all $\alpha \leq \tau$.
\newline
Of course, we begin the induction, by letting $I_0 = \emptyset$
and $M_0=0$. For the inductive step of the proof, we assume that
$\beta$ is an ordinal and the pairs $(I_{\alpha},M_{\alpha})$
have been constructed for all ordinals $\alpha < \beta$, so that
properties (i)-(vii) above hold. In the case where $\beta$ is a
limit ordinal, we define
$I_{\beta} = \bigcup_{\alpha < \beta} I_{\alpha}$ and
$M_{\beta} = \bigcup_{\alpha < \beta} M_{\alpha}$. Then, properties
(i)-(vii) still hold, since
\[ M_{\beta} \otimes_k B =
   {\textstyle{\bigcup_{\alpha < \beta}}} (M_{\alpha} \otimes_k B) =
   {\textstyle{\bigcup_{\alpha < \beta}}} P(I_{\alpha}) =
   P(I_{\beta}) \subseteq P(I) = M \otimes_k B . \]
We now consider the case where $\beta = \alpha +1$ is a successor
ordinal and note that the quotient $M/M_{\alpha}$ is a cofibrant
module, since
\begin{equation}
 (M/M_{\alpha}) \otimes_k B =
 (M \otimes_k B)/(M_{\alpha} \otimes_k B) =
 P(I)/P(I_{\alpha}) = P(I \setminus I_{\alpha})
\end{equation}
is a projective $kG$-module. It follows from the proof of Proposition \ref{prop:KapC},
applied to the cofibrant module $M/M_{\alpha}$, its element
$x_{\alpha}+M_{\alpha}$ and the decomposition of
$(M/M_{\alpha}) \otimes_k B$ displayed in (1) above, that there exist
a subset $I_{\alpha +1} \subseteq I$ containing $I_{\alpha}$ and a
$kG$-submodule $M_{\alpha +1} \subseteq M$ containing $M_{\alpha}$
and $x_{\alpha}$, such that
$\mbox{card} \, (I_{\alpha +1} \setminus I_{\alpha}) \leq \lambda$,
the $kG$-submodule $M_{\alpha +1}/M_{\alpha} \subseteq M/M_{\alpha}$
is $\lambda$-generated and
$(M_{\alpha +1}/M_{\alpha}) \otimes_k B =
 P(I_{\alpha +1} \setminus I_{\alpha})$.
Since $M_{\alpha} \otimes_k B = P(I_{\alpha})$ and
\[ (M_{\alpha +1} \otimes_k B)/(M_{\alpha} \otimes_k B) =
   (M_{\alpha +1}/M_{\alpha}) \otimes_k B =
   P(I_{\alpha +1} \setminus I_{\alpha}) , \]
it is easily seen that
$M_{\alpha +1} \otimes_k B = P(I_{\alpha +1}) \subseteq
 P(I) = M \otimes_k B$.
Indeed, $M_{\alpha +1} \otimes_k B$ is clearly contained
in $P(I_{\alpha +1})$ and the embedding
$\imath : M_{\alpha +1} \otimes_k B \hookrightarrow
          P(I_{\alpha +1})$
fits into the commutative diagram with exact rows
\[
\begin{array}{ccccccccc}
 0 & \longrightarrow & M_{\alpha} \otimes_k B
   & \longrightarrow & M_{\alpha +1} \otimes_k B
   & \longrightarrow
   & (M_{\alpha +1} \otimes_k B)/(M_{\alpha} \otimes_k B)
   & \longrightarrow & 0 \\
   & & \downarrow & & \!\!\! {\scriptstyle{\imath}} \downarrow
   & & \downarrow & & \\
 0 & \longrightarrow & P(I_{\alpha}) & \longrightarrow
   & P(I_{\alpha +1}) & \longrightarrow
   & P(I_{\alpha +1} \setminus I_{\alpha})
   & \longrightarrow & 0
\end{array}
\]
whose end vertical arrows are both bijective. The inductive
step of the construction is therefore complete.

Having constructed the family of pairs
$(I_{\alpha},M_{\alpha})_{\alpha \leq \tau}$ as above, we note that
the $kG$-submodule $M_{\tau} = \bigcup_{\alpha < \tau}M_{\alpha}$
contains all of the generators $(x_{\alpha})_{\alpha < \tau}$ of
$M$ and hence $M_{\tau}=M$. Moreover, for any $\alpha < \tau$ the
$kG$-module
$(M_{\alpha +1}/M_{\alpha}) \otimes_k B =
 P(I_{\alpha +1} \setminus I_{\alpha})$
is projective, so that $M_{\alpha +1}/M_{\alpha}$ is cofibrant.
Hence, $(M_{\alpha})_{\alpha \leq \tau}$ is the required continuous
ascending filtration of $M$ by $\lambda$-generated cofibrant modules.
\end{proof}

\noindent
We can now state and prove the main result of this section.

\begin{Theorem}\label{thm:cp}
The pair $\left( {\tt Cof}(kG),{\tt Cof}(kG)^{\perp} \right)$
is a cotorsion pair in the category of $kG$-modules, which is
cogenerated by a set. In particular, it is a complete cotorsion
pair.
\end{Theorem}

\begin{proof}
In view of Proposition \ref{prop:caf}, there exists a cardinal number
$\lambda$, such that any cofibrant module admits a continuous
ascending filtration by $\lambda$-generated cofibrant modules.
We consider a {\em set} ${\mathcal S}$ of representatives of
the isomorphism classes of all $\lambda$-generated cofibrant
modules and prove that
\begin{equation}
 \left( {\tt Cof}(kG),{\tt Cof}(kG)^{\perp} \right) =
 \left( ^{\perp} \! \left( \mathcal{S}^{\perp} \right) \! ,
 \mathcal{S}^{\perp} \right)
\end{equation}
is the cotorsion pair cogenerated by $\mathcal{S}$. First of
all, we shall prove that ${\tt Cof}(kG)^{\perp} = \mathcal{S}^{\perp}$.
Since $\mathcal{S} \subseteq {\tt Cof}(kG)$, it follows that
${\tt Cof}(kG)^{\perp} \subseteq \mathcal{S}^{\perp}$. In
order to prove the reverse inclusion, we consider a module
$N \in \mathcal{S}^{\perp}$. Then, the functor
${\rm Ext}^1_{kG}(\_\!\_,N)$ vanishes on all modules in
$\mathcal{S}$ and hence on all $\lambda$-generated cofibrant
modules. Since any cofibrant module admits a continuous
ascending filtration by $\lambda$-generated cofibrant modules,
Eklof's lemma \cite[Theorem 7.3.4]{EJ} implies that the functor
${\rm Ext}^1_{kG}(\_\!\_,N)$ vanishes on all cofibrant modules;
hence, $N \in {\tt Cof}(kG)^{\perp}$. We have thus proved the
equality ${\tt Cof}(kG)^{\perp} = \mathcal{S}^{\perp}$. We shall
now prove that
${\tt Cof}(kG) = \! \, ^{\perp} \! \left( \mathcal{S}^{\perp} \right)$.
We note that the equality ${\tt Cof}(kG)^{\perp} = \mathcal{S}^{\perp}$
implies that
\[ {\tt Cof}(kG) \subseteq
   \! \, ^{\perp} \! \left( {\tt Cof}(kG)^{\perp} \right) =
   \! \, ^{\perp} \! \left( \mathcal{S}^{\perp} \right) . \]
Conversely, $\mathcal{S}$ is contained in ${\tt Cof}(kG)$ and
the latter class is closed under direct summands and continuous
ascending filtrations; cf.\ Lemma 2.1. Then,
\cite[Corollary 7.3.5]{EJ} shows that we also have
$^{\perp} \! \left( \mathcal{S}^{\perp} \right) \subseteq
 {\tt Cof}(kG)$.

Having proved (2), it follows from \cite[Theorem 10]{ETz}
that $\left( {\tt Cof}(kG),{\tt Cof}(kG)^{\perp} \right)$
is a complete cotorsion pair.
\end{proof}

\section{Basic properties of cofibrant modules}

\noindent
In this section, we describe a few basic properties of the
class ${\tt Cof}(kG)$ of Benson's cofibrant modules and its
right Ext$^1$-orthogonal ${\tt Cof}(kG)^{\perp}$.

Since all projective $kG$-modules are cofibrant, a standard
argument shows that the following two conditions are equivalent
for a $kG$-module $M$:

(i) There is a non-negative integer $n$ and an exact sequence of
$kG$-modules
\[ 0 \longrightarrow M_n \longrightarrow \cdots \longrightarrow M_1
     \longrightarrow M_0 \longrightarrow N \longrightarrow 0 , \]
such that $M_0,M_1, \ldots ,M_n \in {\tt Cof}(kG)$.

(ii) $\mbox{pd}_{kG}(M \otimes_k B) < \infty$.
\newline
If these conditions are satisfied, we say that the $kG$-module $M$
has finite cofibrant dimension; we denote by $\overline{\tt Cof}(kG)$
the class of these $kG$-modules. It is clear that
${\tt Cof}(kG) \subseteq \overline{\tt Cof}(kG)$. An interesting
criterion for a module in $\overline{\tt Cof}(kG)$ to be cofibrant
was obtained by Benson. For the reader's convenience, we provide a
detailed account of Benson's argument below.

\begin{Lemma}\label{lem:Ben} (cf.\ \cite[Lemma 4.5(ii)]{Ben})
Let $M \in \overline{\tt Cof}(kG)$ and assume that
${\rm Ext}^1_{kG}(M,N)=0$ for any $kG$-module $N$
of finite projective dimension. Then, $M \in {\tt Cof}(kG)$.
\end{Lemma}

\begin{proof}
Assuming that $M$ is not cofibrant, let
$n = \mbox{pd}_{kG}(M \otimes_k B) > 0$ and consider a short
exact sequence of $kG$-modules
\begin{equation}
 0 \longrightarrow K \longrightarrow P \longrightarrow M
   \longrightarrow 0 ,
\end{equation}
where $P$ is projective. Tensoring that exact sequence with
the $k$-free $kG$-module $B$, it follows that
$\mbox{pd}_{kG}(K \otimes_k B) = n-1$. We note that there is
a commutative diagram with exact rows
\[
\begin{array}{ccccccccc}
 0 & \longrightarrow & K \otimes_k \overline{B}
   & \stackrel{\imath}{\longrightarrow}
   & P \otimes_k \overline{B} & \longrightarrow
   & M \otimes_k \overline{B} & \longrightarrow & 0 \\
 & & \!\!\! {\scriptstyle{\varphi}} \downarrow
 & & \!\!\! {\scriptstyle{f}} \downarrow
 & & \parallel & & \\
 0 & \longrightarrow & M
   & \stackrel{\jmath}{\longrightarrow}
   & M \otimes_k B & \longrightarrow
   & M \otimes_k \overline{B} & \longrightarrow & 0
\end{array}
\]
Here, the top row is obtained by tensoring (3) with the
$k$-free $kG$-module $\overline{B}$, the bottom row is
obtained by tensoring with $M$ the $k$-split short exact
sequence of $kG$-modules
\[  0 \longrightarrow k \longrightarrow B
      \longrightarrow \overline{B} \longrightarrow 0 , \]
$f$ is obtained from the projectivity of the $kG$-module
$P \otimes_k \overline{B}$ and $\varphi$ is the restriction
of $f$. We also consider a projective $kG$-module $Q$
admitting a surjective $kG$-linear map
$\pi : Q \longrightarrow M$ and the commutative diagram
with exact rows
\[
\begin{array}{ccccccccc}
 0 & \longrightarrow & Q \oplus \left(
     K \otimes_k \overline{B} \right)
   & \stackrel{1 \oplus \imath}{\longrightarrow}
   & Q \oplus \left( P \otimes_k \overline{B} \right)
   & \longrightarrow & M \otimes_k \overline{B}
   & \longrightarrow & 0 \\
 & & \!\!\! {\scriptstyle{\varphi'}} \downarrow
 & & \!\!\! {\scriptstyle{f'}} \downarrow
 & & \parallel & & \\
 0 & \longrightarrow & M
   & \stackrel{\jmath}{\longrightarrow}
   & M \otimes_k B & \longrightarrow
   & M \otimes_k \overline{B} & \longrightarrow & 0
\end{array}
\]
where $\varphi' = (\pi , \varphi)$ and
$f' = (\jmath \circ \pi , f)$. Since $\varphi'$ is surjective,
an application of the snake lemma shows that $f'$ is also
surjective and $\mbox{ker} \, \varphi' \cong \mbox{ker} \, f'$.
Our assumption that $M \in \overline{\tt Cof}(kG)$ and the
projectivity of $Q \oplus \left( P \otimes_k \overline{B} \right)$
imply that the $kG$-module
$\mbox{ker} \, \varphi' \cong \mbox{ker} \, f'$ has finite
projective dimension. Invoking the additional hypothesis on
$M$, we conclude that the epimorphism $\varphi'$ splits, so
that $M$ is a direct summand of
$Q \oplus \left( K \otimes_k \overline{B} \right)$. Hence,
$M \otimes_k B$ is a direct summand of
$\left[ Q \oplus \left( K \otimes_k \overline{B} \right)
 \right] \otimes_k B = \left( Q \otimes_k B \right) \oplus
 \left( K \otimes_k \overline{B} \otimes_k B \right)$.
This is a contradiction, since $M \otimes_k B$ has projective
dimension equal to $n$, whereas
$\left( Q \otimes_k B \right) \oplus
 \left( K \otimes_k \overline{B} \otimes_k B \right)$
has projective dimension $\leq n-1$.
\end{proof}

\noindent
Cofibrant modules are closely related to Gorenstein projective
modules. Recall that an acyclic complex of projective $kG$-modules
is called totally acyclic if it remains acyclic after applying
the functor ${\rm Hom}_{kG}(\_\!\_,P)$ for any projective
$kG$-module $P$. A $kG$-module $M$ is called Gorenstein projective
if it is a cokernel of a totally acyclic complex of projective
$kG$-modules. Let ${\tt GProj}(kG)$ denote the class of Gorenstein
projective $kG$-modules. We refer to \cite{EJ, Hol} for this notion
(which may be defined over any ring). Cornick and Kropholler showed
in \cite[Theorem 3.5]{CK1} that any cofibrant module is Gorenstein
projective, so that we always have an inclusion
${\tt Cof}(kG) \subseteq {\tt GProj}(kG)$.

\begin{Proposition}\label{prop:cof}
(i) ${\tt Cof}(kG) = \overline{\tt Cof}(kG) \cap {\tt GProj}(kG)$.

(ii) The class ${\tt Cof}(kG)$ is closed under extensions and
kernels of epimorphisms.

(iii) Any cofibrant module is a cokernel of an acyclic complex
of projective $kG$-modules, all of whose cokernels are cofibrant.
\end{Proposition}

\begin{proof}
(i) Since ${\tt Cof}(kG) \subseteq \overline{\tt Cof}(kG)$ and
${\tt Cof}(kG) \subseteq {\tt GProj}(kG)$, it is clear that
${\tt Cof}(kG) \subseteq \overline{\tt Cof}(kG) \cap {\tt GProj}(kG)$.
Conversely, if $M \in \overline{\tt Cof}(kG) \cap {\tt GProj}(kG)$,
then ${\rm Ext}^1_{kG}(M,N)=0$ for any $kG$-module $N$ of finite
projective dimension; cf.\ \cite[Theorem 2.20]{Hol}. Hence, Lemma \ref{lem:Ben}
implies that $M$ is cofibrant, as needed.

(ii) We consider a short exact sequence of $kG$-modules
\[ 0 \longrightarrow M' \longrightarrow M \longrightarrow M''
     \longrightarrow 0 \]
and assume that $M''$ is cofibrant. Then, the $kG$-module
$M'' \otimes_k B$ is projective and hence the associated short
exact sequence of $kG$-modules
\[ 0 \longrightarrow M' \otimes_k B \longrightarrow M \otimes_k B
     \longrightarrow M'' \otimes_k B \longrightarrow 0 \]
splits. We conclude that the $kG$-module $M' \otimes_k B$ is
projective if and only if the $kG$-module $M \otimes_k B$ is
projective. In other words, $M'$ is cofibrant if and only if
$M$ is cofibrant.

(iii) Let $M$ be a cofibrant module. Since all projective
$kG$-modules are cofibrant, it follows from (ii) above that
the cokernels of any projective resolution of $M$ are cofibrant.
This provides us with the left half of the required acyclic
complex of projective $kG$-modules. On the other hand, being
cofibrant, the $kG$-module $M$ is Gorenstein projective. Hence,
there exists a short exact sequence of $kG$-modules
\[ 0 \longrightarrow M \longrightarrow P \longrightarrow N
     \longrightarrow 0 , \]
where $P$ is projective (and hence cofibrant) and $N$ is
Gorenstein projective. It follows readily that
$N \in \overline{\tt Cof}(kG)$, so that (i) above implies
that the $kG$-module $N$ is actually cofibrant. We may repeat
the same argument with $N$ in place of $M$ and construct
inductively the right half of the required acyclic complex
of projective $kG$-modules, by splicing short exact sequences
with cofibrant end terms and projective middle term.
\end{proof}

\begin{Proposition}\label{prop:Proj}
The cotorsion pair of Theorem 2.4 is hereditary and projective.
\end{Proposition}

\begin{proof}
The cotorsion pair is hereditary, since its left-hand side
class ${\tt Cof}(kG)$ is closed under kernels of epimorphisms;
cf.\ Proposition 3.2(ii).

We shall now prove that the cotorsion pair is projective, i.e.\
that ${\tt Cof}(kG)\cap{\tt Cof}(kG)^\perp = {\tt Proj}(kG)$.
Since all cofibrant modules are Gorenstein projective, it is
clear that ${\tt Proj}(kG) \subseteq {\tt Cof}(kG)^\perp$. On
the other hand, any projective $kG$-module is cofibrant and hence
${\tt Proj}(kG)\subseteq {\tt Cof}(kG)\cap{\tt Cof}(kG)^\perp$.
For the reverse inclusion, let $M$ be a $kG$-module contained in
the intersection ${\tt Cof}(kG)\cap{\tt Cof}(kG)^\perp$. Since
$M$ is cofibrant, Proposition 3.2(iii) implies the existence of
a short exact sequence of $kG$-modules
\[ 0 \longrightarrow M \longrightarrow P \longrightarrow N
     \longrightarrow 0 , \]
where $P$ is projective and $N$ is cofibrant. Our hypothesis
that $M \in {\tt Cof}(kG)^{\perp}$ implies that the short exact
sequence above splits, so that $M$ is a direct summand of the
projective $kG$-module $P$. In particular, the $kG$-module $M$
is projective, as needed.
\end{proof}

\noindent
The following simple property regarding the behaviour of the class
of cofibrant modules with respect to tensor products will be useful
in the sequel.

\begin{Lemma}\label{lem:obser}
Let $M,N$ be two $kG$-modules.

(i) If $M$ is cofibrant, then $M$ is $k$-projective.

(ii) If $M$ is $k$-projective and $N$ is cofibrant,
then the diagonal $kG$-module $M\otimes_k N$ is cofibrant.

(iii) If $M$ and $N$ are cofibrant, then the diagonal
$kG$-module $M\otimes_k N$ is cofibrant.
\end{Lemma}

\begin{proof}
(i) Since $k$ is a direct summand of $B$ as a a $k$-module, it follows
that $M$ is a direct summand of $M \otimes_kB$ as a $k$-module. Since
$M \otimes_kB$ is projective (as a $kG$-module and hence) as a $k$-module,
we conclude that $M$ is also projective as a $k$-module.

(ii) The assumption on $M$ implies that for any projective $kG$-module
$P$ the diagonal $kG$-module $M \otimes_k P$ is projective as well. In
particular, $(M \otimes_k N) \otimes_k B = M \otimes_k (N \otimes_k B)$
is a projective $kG$-module and hence $M\otimes_k N$ is cofibrant.

(iii) This is an immediate consequence of (i) and (ii).
\end{proof}

\noindent
The orthogonal class ${\tt Cof}(kG)^{\perp}$ is of some interest,
as it determines the class of cofibrant modules as its left
Ext$^1$-orthogonal. Here are some basic properties of that class.

\begin{Proposition}\label{prop:cof+}
The class ${\tt Cof}(kG)^{\perp}$ is thick, i.e.\ it is closed
under direct summands and has the 2-out-of-3 property for short
exact sequences. It contains all $kG$-modules of finite flat or
injective dimension and all (diagonal) $kG$-modules of the form
${\rm Hom}_k(B,N)$, where $N$ is a $kG$-module.
\end{Proposition}

\begin{proof}
The class ${\tt Cof}(kG)^{\perp}$ is clearly closed under direct
summands. Since the cotorsion pair
$\left( {\tt Cof}(kG),{\tt Cof}(kG)^{\perp} \right)$ is hereditary,
${\tt Cof}(kG)^{\perp}$ is closed under extensions and cokernels of
monomorphisms. We prove that it is also closed under kernels of
epimorphisms. Let
\begin{equation}
 0 \longrightarrow N' \longrightarrow N \longrightarrow N''
   \longrightarrow 0
\end{equation}
be a short exact sequence of $kG$-modules and assume that
$N,N'' \in {\tt Cof}(kG)^{\perp}$. In order to show that
$N' \in {\tt Cof}(kG)^{\perp}$, we consider a cofibrant module
$M$ and invoke Proposition \ref{prop:cof}(iii) to find a short exact
sequence of $kG$-modules
\[ 0 \longrightarrow M \longrightarrow P \longrightarrow M'
     \longrightarrow 0 , \]
where $P$ is projective and $M'$ is cofibrant. Then,
${\rm Ext}^1_{kG}(M,N') \cong {\rm Ext}^2_{kG}(M',N')$ and hence
it suffices to prove that the abelian group ${\rm Ext}^2_{kG}(M',N')$
is trivial. This follows from the long exact sequence associated
with (4)
\[ \cdots \longrightarrow {\rm Ext}^1_{kG}(M',N'')
   \longrightarrow {\rm Ext}^2_{kG}(M',N')
   \longrightarrow {\rm Ext}^2_{kG}(M',N) \longrightarrow \cdots ,\]
since both groups ${\rm Ext}^1_{kG}(M',N'')$ and
${\rm Ext}^2_{kG}(M',N)$ are trivial.

We note that ${\tt Cof}(kG)^{\perp}$ contains the class of flat
modules; this follows from \cite[Proposition 8.2]{Ste} and
\cite[Theorem 4.4]{SS}. Of course, ${\tt Cof}(kG)^{\perp}$
also contains the class of injective modules. Using induction
and the 2-out-of-3 property for short exact sequences established
above, it follows that ${\tt Cof}(kG)^{\perp}$ contains all
$kG$-modules of finite flat dimension and all $kG$-modules of
finite injective dimension. Finally, if $N$ is any $kG$-module,
the Hom-tensor adjunction for diagonal $kG$-modules implies that
the functor
\[ {\rm Ext}^1_{kG}(\_\!\_ \otimes_k B,N) \cong
   {\rm Ext}^1_{kG}(\_\!\_,{\rm Hom}_k(B,N)) \]
vanishes on all cofibrant modules and hence
${\rm Hom}_k(B,N) \in {\tt Cof}(kG)^{\perp}$.
\end{proof}

\section{Cofibrant and Gorenstein projective modules}

\noindent
In this section, we elaborate on the relation between Benson's
cofibrant and Gorenstein projective modules. Our goal is to
enlarge (in principle, at least) the class of groups over which
these classes are known to coincide.

Let us denote by $\mathcal{P}(kG)$ the class of those
$kG$-modules that appear as cokernels of acyclic complexes
of projective $kG$-modules. It is clear that
${\tt GProj}(kG) \subseteq \mathcal{P}(kG)$ and equality holds
if and only if any acyclic complex of projective $kG$-modules
is totally acyclic. It is conjectured in \cite{DT} that
$\mathcal{P}(\mathbb{Z}G) \subseteq {\tt Cof}(\mathbb{Z}G)$,
so that
${\tt Cof}(\mathbb{Z}G) = {\tt GProj}(\mathbb{Z}G) =
 \mathcal{P}(\mathbb{Z}G)$.
If these equalities hold, then the study of Gorenstein
projective $\mathbb{Z}G$-modules is reduced to questions
involving classical homological algebra notions. It is
proved in [loc.cit., Corollary C] that the equalities
above do hold for any group in Kropholler's class
${\scriptstyle{{\bf LH}}}\mathfrak{F}$. The question as to
whether Gorenstein projective modules are necessarily cofibrant
has been studied over more general coefficient rings. For
example, it is proved in \cite[Corollary 5.5]{Bi} that
$\mathcal{P}(kG) \subseteq {\tt Cof}(kG)$, if the ring
$k$ has finite global dimension and $G$ is either an
${\scriptstyle{{\bf LH}}}\mathfrak{F}$-group or a group of
type $\Phi$.

The following remark shows that certain restrictions on the
coefficient ring $k$ are necessary for the equalities
${\tt Cof}(kG) = {\tt GProj}(kG) = \mathcal{P}(kG)$
to hold.

\vspace{0.1in}

\noindent
{\bf Remarks 4.1}
(i) If $G$ is a group and ${\tt Cof}(kG) = \mathcal{P}(kG)$,
then any acyclic complex of projective $k$-modules is contractible.
Indeed, let $P_*$ be an acyclic complex of projective $k$-modules
with cokernels $(C_nP)_n$. Then, $\mbox{ind}_1^GP_*$ is an acyclic
complex of projective $kG$-modules with cokernels
$(\mbox{ind}_1^GC_nP)_n$. Since ${\tt Cof}(kG) = \mathcal{P}(kG)$,
these cokernels are cofibrant $kG$-modules and hence the restricted
$k$-module $\mbox{res}_1^G\mbox{ind}_1^GC_nP$ is projective for all
$n$; cf.\ Lemma 3.4(i). It follows that $C_nP$, being a direct
summand of $\mbox{res}_1^G\mbox{ind}_1^GC_nP$, is a projective
$k$-module for all $n$. Therefore, the complex $P_*$ is indeed
contractible.

(ii) If $G$ is a group and ${\tt Cof}(kG) = {\tt GProj}(kG)$,
then any totally acyclic complex of projective $k$-modules is
contractible. This follows using exactly the same argument as
in (i) above, since the induction functor preserves the class
of totally acyclic complexes of projective modules.
\addtocounter{Lemma}{1}

(iii) For an explicit example where the inclusion
${\tt Cof}(kG) \subseteq {\tt GProj}(kG)$ is proper, let
$k = R[x]/(x^2)$ be the ring of dual numbers over a commutative
ring $R$. Then, the $k$-module $N = (x)/(x^2)$ is Gorenstein
projective, but not projective. Therefore, for any group $G$,
the induced $kG$-module $M = \mbox{ind}_1^GN$ is Gorenstein
projective, but not cofibrant.

\vspace{0.1in}

\noindent
In order to find examples of commutative rings $k$ over which
all acyclic complexes of projective modules are contractible
(so that one may hope, in view of Remark 4.1, that there is an
equality ${\tt Cof}(kG) = {\tt GProj}(kG)$ for certain groups
$G$), we recall that a $k$-module is pure-projective if it is
a direct summand of a suitable direct sum of finitely presented
$k$-modules. Then, the following two conditions are easily seen
to be equivalent:

(i) All pure-projective $k$-modules have finite projective
dimension.

(ii) All finitely presented $k$-modules have finite projective
dimension and there is a uniform upper bound on their projective
dimension.

\begin{Lemma}
Assume that $k$ satisfies one of the following two conditions:

(i) $k$ has finite weak global dimension or

(ii) all pure-projective $k$-modules have finite projective
dimension.
\newline
Then, any acyclic complex of projective $k$-modules is contractible.
\end{Lemma}

\begin{proof}
The conclusion under hypothesis (i) follows as in
\cite[Proposition 1.1]{ET1}, from the result by Benson and
Goodearl \cite[Theorem 2.5]{BG} and Neeman \cite[Remark 2.15]{N2}
on the contractibility of pure acyclic complexes of projective
modules; see also \cite[Proposition 7.6]{CH}. The conclusion
under hypothesis (ii) follows by applying \cite[Remark 2.11(i)]{ET2}
to the cotorsion pair $({\tt Proj}(k),k\mbox{-Mod})$.
\end{proof}

\noindent
{\bf Remark 4.3}
If $k$ is $\aleph_n$-Noetherian for some integer $n \geq 0$,
the hypothesis in Lemma 4.2(ii) that all pure-projective
$k$-modules have finite projective dimension implies that the
global dimension of $k$ is finite (so that this hypothesis is
strictly stronger than that in Lemma 4.2(i)). To verify this
claim, fix a positive integer $s$, such that $\mbox{pd}_kC \leq s$
for any finitely presented $k$-module $C$. In particular, it
follows that $\mbox{pd}_kI \leq s-1$ for any finitely generated
ideal $I \subseteq k$. Using the $\aleph_n$-Noetherian hypothesis
on $k$ and invoking \cite[Corollary 2.7]{ET2}, we conclude that
$\mbox{id}_kN \leq s+n+1$ for any $k$-module $N$ and hence
$\mbox{gl.dim} \, k \leq s+n+1$.
\addtocounter{Lemma}{1}

\vspace{0.1in}

\noindent
Having fixed the commutative ring $k$, let
$\mathfrak{C} = \mathfrak{C}(k)$ be the class consisting
of those groups $H$ that have the following property:
{\em If $G$ is a group containing $H$ as a subgroup
and $M \in \mathcal{P}(kG)$, then the $kH$-module
$\mbox{res}_H^G(M \otimes_k B(G)) =
 \mbox{res}_H^GM \otimes_k \mbox{res}_H^GB(G)$
is projective.} In particular, letting $G=H$, we conclude
that for any $\mathfrak{C}$-group $H$ we have an equality
$\mathcal{P}(kH) = {\tt Cof}(kH)$. Here are some basic
properties of $\mathfrak{C}$.

\begin{Lemma}
The following conditions are equivalent for a commutative
ring $k$:

(i) All finite groups are contained in $\mathfrak{C}$.

(ii) The class $\mathfrak{C}$ is non-empty.

(iii) Any acyclic complex of projective $k$-modules is
contractible.
\end{Lemma}

\begin{proof}
It is clear that (i)$\rightarrow$(ii). In order to show
that (ii)$\rightarrow$(iii), assume that $H$ is a
$\mathfrak{C}$-group. Then, as we have noted above, there
is an equality $\mathcal{P}(kH) = {\tt Cof}(kH)$; hence,
Remark 4.1(i) implies that (iii) holds. Finally, assume
that (iii) holds and let $F$ be a finite group. Let $G$
be a group containing $F$ as a subgroup and consider a
$kG$-module $M \in \mathcal{P}(kG)$. Our assumption on
$k$ implies that $M$ is projective as a $k$-module, whereas
$B(G)$ is projective as a $kF$-module; cf.\ \cite{KT}.
It follows that $M \otimes_k B(G)$ is projective as a
$kF$-module. Since this is the case for any group $G$
containing $F$ and any $kG$-module $M \in \mathcal{P}(kG)$,
it follows that $F$ is a $\mathfrak{C}$-group.
\end{proof}

\begin{Lemma}
The class $\mathfrak{C}$ is subgroup-closed.
\end{Lemma}

\begin{proof}
Let $F$ be a subgroup of a $\mathfrak{C}$-group $H$. In order
to show that $F \in \mathfrak{C}$, we consider a group $G$
containing $F$ and a module $M \in \mathcal{P}(kG)$. Let
$\Gamma = G*_FH$ be the amalgamated free product of $G$
and $H$ along their common subgroup $F$ and consider the
$k\Gamma$-module
$N = \mbox{ind}_G^{\Gamma}M \in \mathcal{P}(k\Gamma)$. Since
$H$ is a $\mathfrak{C}$-group, the $kH$-module
$\mbox{res}_H^{\Gamma}N \otimes_k \mbox{res}_H^{\Gamma}B(\Gamma)$
is projective. Restricting further to $F$, we conclude that
the $kF$-module
$\mbox{res}_F^{\Gamma}N \otimes_k \mbox{res}_F^{\Gamma}B(\Gamma)$
is projective as well. We note that the $kG$-module $M$ is a
direct summand of
$\mbox{res}_G^{\Gamma}\mbox{ind}_G^{\Gamma}M =
 \mbox{res}_G^{\Gamma}N$
and hence the $kF$-module $\mbox{res}_F^GM$ is a direct summand
of $\mbox{res}_F^{\Gamma}N$. On the other hand, $B(G)$ is a
direct summand of $\mbox{res}_G^{\Gamma}B(\Gamma)$ and hence
$\mbox{res}_F^GB(G)$ is a direct summand of
$\mbox{res}_F^{\Gamma}B(\Gamma)$.
It follows that $\mbox{res}_F^GM \otimes_k \mbox{res}_F^GB(G)$
is a direct summand of
$\mbox{res}_F^{\Gamma}N \otimes_k \mbox{res}_F^{\Gamma}B(\Gamma)$,
so that $\mbox{res}_F^GM \otimes_k \mbox{res}_F^GB(G)$ is a
projective $kF$-module as well.
\end{proof}

\noindent
The following result will be useful, in order to exhibit
a plethora of $\mathfrak{C}$-groups.

\begin{Theorem}
The class $\mathfrak{C}$ is closed under the operation
${\scriptstyle{{\bf LH}}} \cup \Phi \cup \Phi_{flat}$.
\end{Theorem}

\begin{proof}
We have to show that
$\mathfrak{C} = {\scriptstyle{{\bf LH}}}\mathfrak{C}
 \cup \Phi\mathfrak{C} \cup \Phi_{flat}\mathfrak{C}$,
i.e.\ that
${\scriptstyle{{\bf LH}}}\mathfrak{C} \cup
 \Phi\mathfrak{C} \cup \Phi_{flat}\mathfrak{C}
 \subseteq \mathfrak{C}$.
Let $H$ be a group contained in
${\scriptstyle{{\bf LH}}}\mathfrak{C} \cup
 \Phi\mathfrak{C} \cup \Phi_{flat}\mathfrak{C}$.
In order to prove that $H$ is a $\mathfrak{C}$-group,
we consider a group $G$ with $H \subseteq G$ and
a $kG$-module $M \in \mathcal{P}(kG)$. We have to show
that the $kH$-module $\mbox{res}_H^G(M \otimes_k B(G))$
is projective. This is a consequence of the following
result.
\end{proof}

\begin{Proposition}
Let $G$ be a group. Then, for any $kG$-module
$M \in \mathcal{P}(kG)$ and any subgroup
$H \subseteq G$ with
$H \in {\scriptstyle{{\bf LH}}}\mathfrak{C} \cup
 \Phi\mathfrak{C} \cup \Phi_{flat}\mathfrak{C}$,
the $kH$-module $\mbox{res}_H^G(M \otimes_k B(G))$
is projective.
\end{Proposition}

\begin{proof}
We let $N = M \otimes_k B(G)$ and distinguishing four
cases, according to whether (i)
$H \in {\scriptstyle{{\bf H}}}\mathfrak{C}$, (ii)
$H \in {\scriptstyle{{\bf LH}}}\mathfrak{C}$, (iii)
$H \in \Phi\mathfrak{C}$ and (iv) $H \in \Phi_{flat}\mathfrak{C}$.

(i) Assume that $H \in {\scriptstyle{{\bf H}}}\mathfrak{C}$
and use induction on the ordinal $\alpha$, which is such that
$H \in {\scriptstyle{{\bf H}}}_{\alpha}\mathfrak{C}$. If
$\alpha = 0$, then $H \in \mathfrak{C}$ and hence
$\mbox{res}_H^GN \in {\tt Proj}(kH)$, in view of the very
definition of $\mathfrak{C}$-groups. Let $\alpha >0$ and
assume that the result is known for all
${\scriptstyle{{\bf H}}}_{\beta}\mathfrak{C}$-subgroups
$K \subseteq G$ and all ordinals $\beta$ with $\beta < \alpha$.
The group $H$ acts cellularly on a finite dimensional
contractible CW-complex $X$, with each cell stabilizer
contained in ${\scriptstyle{{\bf H}}}_{\beta}\mathfrak{C}$
for a suitable $\beta < \alpha$. Then, there is an induced
short exact sequence of $kH$-modules
\[ 0 \longrightarrow N_d \longrightarrow \cdots
     \longrightarrow N_1 \longrightarrow N_0
     \longrightarrow \mbox{res}_H^GN \longrightarrow 0 , \]
where $d$ is the dimension of $X$ and each $N_i$ is a direct
sum of modules of the form $\mbox{ind}_K^H\mbox{res}_K^GN$,
for a suitable
${\scriptstyle{{\bf H}}}_{\beta}\mathfrak{C}$-subgroup
$K \subseteq H$ for some $\beta < \alpha$. In view of the
induction hypothesis, we have $N_i \in {\tt Proj}(kH)$ for
all $i=0,1, \ldots ,d$ and hence
$\mbox{pd}_{kH} \mbox{res}_H^GN \leq d$. This inequality
holds for any $M \in \mathcal{P}(kG)$. Since the $kG$-module
$M$ is the $d$-th syzygy in a projective resolution of
another $kG$-module $M' \in \mathcal{P}(kG)$, a standard
argument shows that the $kH$-module $\mbox{res}_H^GN$ is
actually projective.

(ii) Assume that $H \in {\scriptstyle{{\bf LH}}}\mathfrak{C}$
and consider a finitely generated subgroup $K \subseteq H$.
Then, there exists an
${\scriptstyle{{\bf H}}}\mathfrak{C}$-subgroup $K' \subseteq H$
with $K \subseteq K'$. Having taken care of all
${\scriptstyle{{\bf H}}}\mathfrak{C}$-subgroups of $H$ in
(i) above, we conclude that
$\mbox{res}_{K'}^GN \in {\tt Proj}(kK')$. Restricting further
to $K$, it follows that $\mbox{res}_K^GN \in {\tt Proj}(kK)$.
In particular, the $kK$-module $\mbox{res}_K^GN$ is flat and
hence the induced $kH$-module $\mbox{ind}_K^H\mbox{res}_K^GN$
is flat as well. Since the $kH$-module $\mbox{res}_H^GN$ is the
filtered colimit of the system $(\mbox{ind}_K^H\mbox{res}_K^GN)_K$,
where $K$ runs through the family of all finitely generated
subgroups of $H$, we conclude that $\mbox{res}_H^GN$ is flat.
This is the case for any $M \in \mathcal{P}(kG)$. Since $M$
is a cokernel of an acyclic complex of projective $kG$-modules
with all cokernels in $\mathcal{P}(kG)$, it follows that
$\mbox{res}_H^GN$ is a cokernel of an acyclic complex of
projective $kH$-modules, all of whose cokernels are flat.
Therefore, the result by Benson, Goodearl and Neeman on the
contractibility of pure acyclic complexes of projective modules
(cf.\ \cite[Theorem 2.5]{BG} and \cite[Remark 2.15]{N2}) implies
that $\mbox{res}_H^GN$ is actually a projective $kH$-module.

(iii) Assume that $H \in \Phi\mathfrak{C}$. For any
$\mathfrak{C}$-subgroup $K \subseteq H$ we know that
$\mbox{res}_K^GN \in {\tt Proj}(kK)$, in view of the
definition of $\mathfrak{C}$-groups; hence,
$\mbox{pd}_{kH} \mbox{res}_H^GN < \infty$. Since this
is the case for any $M \in \mathcal{P}(kG)$ and the class
$\mathcal{P}(kG)$ is closed under direct sums, a standard
argument shows that there is a uniform upper bound on
the projective dimension of $\mbox{res}_H^GN$ for all
$M \in \mathcal{P}(kG)$. We may now conclude, as in the
last part of the proof of (i) above, that the $kH$-module
$\mbox{res}_H^GN$ is actually projective.

(iv) Finally, assume that $H \in \Phi_{flat}\mathfrak{C}$. For
any $\mathfrak{C}$-subgroup $K \subseteq H$ the $kH$-module
$\mbox{res}_K^GN$ is projective and hence flat; therefore,
$\mbox{fd}_{kH} \mbox{res}_H^GN < \infty$. Using the same
argument as in the proof of (iii) above, we conclude that
the $kH$-module $\mbox{res}_H^GN$ is flat. Furthermore,
invoking the result by Benson, Goodearl and Neeman on the
contractibility of pure acyclic complexes of projective modules
(cf.\ \cite[Theorem 2.5]{BG} and \cite[Remark 2.15]{N2}),
as in the  last part of the proof of (ii) above, it follows
that the $kH$-module $\mbox{res}_H^GN$ is actually projective.
\end{proof}

\noindent
Let $\overline{\mathfrak{F}}$ be the closure of the class
$\mathfrak{F}$ of finite groups under the operation
${\scriptstyle{{\bf LH}}} \cup \Phi \cup \Phi_{flat}$. As
we explained in $\S $1.II, groups in $\overline{\mathfrak{F}}$
can be described hierarchically starting from the class of
finite groups.

\begin{Corollary}
Assume that any acyclic complex of projective $k$-modules is
contractible and let $G$ be an $\overline{\mathfrak{F}}$-group.
Then, ${\tt Cof}(kG) = {\tt GProj}(kG) = \mathcal{P}(kG)$.
\end{Corollary}

\begin{proof}
In view of Lemma 4.4, our assumption on $k$ implies that
$\mathfrak{F} \subseteq \mathfrak{C}$. Using Theorem 4.6,
we conclude that
$\overline{\mathfrak{F}} \subseteq \mathfrak{C}$ and hence
$G \in \mathfrak{C}$. The equality
$\mathcal{P}(kG) = {\tt Cof}(kG)$ follows readily, as we
have already noted before Lemma 4.4.
\end{proof}

\begin{Corollary}
Assume that $k$ satisfies one of the following two conditions:

(i) $k$ has finite weak global dimension or

(ii) all pure-projective $k$-modules have finite projective
dimension.
\newline
If $G$ is an $\overline{\mathfrak{F}}$-group, then
${\tt Cof}(kG) = {\tt GProj}(kG) = \mathcal{P}(kG)$.
\end{Corollary}

\begin{proof}
This is an immediate consequence of Corollary 4.8, in view
of Lemma 4.2.
\end{proof}

\noindent
{\bf Remark 4.10.}
Assuming that ${\tt Cof}(kG) = {\tt GProj}(kG)$, the
assertions stated in Proposition \ref{prop:cof}(ii) represent
well-known properties of Gorenstein projective modules, whereas
Proposition \ref{prop:cof}(iii) is essentially a consequence of
the very definition of Gorenstein projective modules. Analogously,
the equality ${\tt Cof}(kG) = {\tt GProj}(kG)$ reduces Lemma
\ref{lem:Ben} to the assertion that, within the class of $kG$-modules
of finite Gorenstein projective dimension, the Gorenstein projective
$kG$-modules are precisely those $kG$-modules which are left
Ext$^1$-orthogonal to the $kG$-modules of finite projective
dimension; see \cite[Theorem 2.10]{Hol}.

\section{Cofibrant modules and induction from subgroups}

\noindent
The restriction functor from a group to a subgroup always maps
cofibrant modules to cofibrant modules. This is clear, since for
any subgroup $H$ of a group $G$ the $kH$-module $B(H)$ is a direct
summand of $\mbox{res}_H^GB(G)$. In this section, we shall examine
the extent to which the induction functor from a subgroup preserves
the class of cofibrant modules.

Having fixed the commutative ring $k$, for any group $G$ we denote
by $\mathcal{P}_0(kG)$ the class of those $kG$-modules that appear
as cokernels of acyclic complexes of projective $kG$-modules, that
are contractible as complexes of $k$-modules. We note that Proposition
3.2(iii) and Lemma 3.4(i) imply that
${\tt Cof}(kG) \subseteq \mathcal{P}_0(kG)$. We consider the group
class $\mathfrak{I} = \mathfrak{I}(k)$, which consists of those
groups $H$ that have the following property: {\em If $G$ is a group
containing $H$ as a subgroup and $M \in \mathcal{P}_0(kH)$, then
the diagonal $kH$-module $M \otimes_k \mbox{res}_H^GB(G)$ is projective.}
The relevance of $\mathfrak{I}$-groups to the study of the behaviour
of cofibrant modules under induction stems from Proposition 5.2
below. The argument is based on the following well-known result,
which is really a consequence of the Hopficity of group algebras.

\begin{Lemma} (cf.\ \cite[Lemma 1.1]{Sw})
Let $H$ be a subgroup of a group $G$ and consider a $kH$-module
$M$ and a $kG$-module $N$. Then, there is an isomorphism of
$kG$-modules between the diagonal $kG$-module
$\mbox{ind}_H^GM \otimes_k N$ and the $kG$-module
$\mbox{ind}_H^G(M \otimes_k \mbox{res}_H^GN)$, which is induced
from the diagonal $kH$-module $M \otimes_k \mbox{res}_H^GN$.
\end{Lemma}

\begin{Proposition}
Let $H$ be an $\mathfrak{I}$-subgroup of a group $G$. Then:

(i) the induction functor $\mbox{ind}_H^G$ maps ${\tt Cof}(kH)$
into ${\tt Cof}(kG)$ and

(ii) the restriction functor $\mbox{res}_H^G$ maps
${\tt Cof}(kG)^{\perp}$ into ${\tt Cof}(kH)^{\perp}$.
\end{Proposition}

\begin{proof}
(i) Let $M$ be a cofibrant $kH$-module; then, as we have noted
above, $M$ is contained in $\mathcal{P}_0(kH)$. Since
$H \in \mathfrak{I}$, the diagonal $kH$-module
$M \otimes_k \mbox{res}_H^GB(G)$ is projective and hence the
induced $kG$-module $\mbox{ind}_H^G(M \otimes_k \mbox{res}_H^GB(G))$
is projective as well. In view of Lemma 5.1, that $kG$-module
is isomorphic with the diagonal $kG$-module
$\mbox{ind}_H^GM \otimes_k B(G)$. It follows that the latter
$kG$-module is projective and hence $\mbox{ind}_H^GM \in {\tt Cof}(kG)$.

(ii) This is an immediate consequence of (i), in view of the
Eckmann-Shapiro isomorphism
$\mbox{Ext}^1_{kG} \! \left( \mbox{ind}_H^G \_\!\_ \, ,
 \_\!\_ \right) \simeq \mbox{Ext}^1_{kH} \! \left( \_\!\_
 \, , \mbox{res}_H^G \_\!\_ \right)$.
\end{proof}

\noindent
The definition of $\mathfrak{I}$-groups is (formally) similar
to the definition of $\mathfrak{C}$-groups given in Section 4;
it is easily seen that $\mathfrak{C} \subseteq \mathfrak{I}$.
We proceed as in Section 4 and begin with some basic properties
of the group class $\mathfrak{I}$.

\begin{Lemma}
The class $\mathfrak{I}$ has the following properties:

(i) It contains all finite groups and

(ii) it is subgroup-closed.
\end{Lemma}

\begin{proof}
(i) Let $H$ be a finite group. In order to show that
$H \in \mathfrak{I}$, consider a group $G$ containing
$H$ as a subgroup and a module $M \in \mathcal{P}_0(kH)$.
Then, $M$ is $k$-projective. Since $B(G)$ is free as a
$kH$-module (cf.\ \cite{KT}), the diagonal $kH$-module
$M \otimes_k \mbox{res}_H^GB(G)$ is projective.

(ii) Let $F$ be a subgroup of an $\mathfrak{I}$-group $H$. In
order to show that $F \in \mathfrak{I}$, we consider a group
$G$ containing $F$ and a module $M \in \mathcal{P}_0(kF)$.
Let $\Gamma = G*_FH$ be the amalgamated free product of $G$
and $H$ along their common subgroup $F$ and consider the induced
$kH$-module $N = \mbox{ind}_F^HM \in \mathcal{P}_0(kH)$.
Since $H$ is an $\mathfrak{I}$-group, the $kH$-module
$N \otimes_k \mbox{res}_H^{\Gamma}B(\Gamma)$ is projective.
Restricting to $F$, we conclude that the $kF$-module
$\mbox{res}_F^HN \otimes_k \mbox{res}_F^{\Gamma}B(\Gamma)$
is projective as well. We note that the $kH$-module $M$
is a direct summand of
$\mbox{res}_F^H\mbox{ind}_F^HM = \mbox{res}_F^HN$. On the
other hand, the $kG$-module $B(G)$ is a direct summand of
$\mbox{res}_G^{\Gamma}B(\Gamma)$ and hence the restricted
$kF$-module $\mbox{res}_F^GB(G)$ is a direct summand of
$\mbox{res}_F^{\Gamma}B(\Gamma)$. It follows that
$M \otimes_k \mbox{res}_F^GB(G)$ is a direct summand of
$\mbox{res}_F^HN \otimes_k \mbox{res}_F^{\Gamma}B(\Gamma)$,
so that $M \otimes_k \mbox{res}_F^GB(G)$ is a projective
$kF$-module as well.
\end{proof}

\noindent
The following result will be used, in order to show that
there is an abundance of $\mathfrak{I}$-groups.

\begin{Theorem}
The class $\mathfrak{I}$ is closed under the operation
${\scriptstyle{{\bf LH}}} \cup \Phi \cup \Phi_{flat}$.
\end{Theorem}

\begin{proof}
We have to show that
$\mathfrak{I} = {\scriptstyle{{\bf LH}}}\mathfrak{I}
 \cup \Phi\mathfrak{I} \cup \Phi_{flat}\mathfrak{I}$,
i.e.\ that
${\scriptstyle{{\bf LH}}}\mathfrak{I} \cup
 \Phi\mathfrak{I} \cup \Phi_{flat}\mathfrak{I}
 \subseteq \mathfrak{I}$.
Let $H$ be a group contained in
${\scriptstyle{{\bf LH}}}\mathfrak{I} \cup
 \Phi\mathfrak{I} \cup \Phi_{flat}\mathfrak{I}$.
In order to prove that $H$ is an $\mathfrak{I}$-group,
we consider a group $G$ with $H \subseteq G$ and
a $kH$-module $M \in \mathcal{P}_0(kH)$. We have to show
that the diagonal $kH$-module $M \otimes_k \mbox{res}_H^GB(G)$
is projective. This is a consequence of the following
result.
\end{proof}

\begin{Proposition}
Let $G$ be a group and $H \subseteq G$ a subgroup
contained in
${\scriptstyle{{\bf LH}}}\mathfrak{I} \cup
 \Phi\mathfrak{I} \cup \Phi_{flat}\mathfrak{I}$.
Then, for any $kH$-module $M \in \mathcal{P}_0(kH)$
the diagonal $kH$-module $M \otimes_k \mbox{res}_H^GB(G)$
is projective.
\end{Proposition}

\begin{proof}
Let $N = M \otimes_k \mbox{res}_H^GB(G)$. We shall proceed
by distinguishing four cases, according to whether (i)
$H \in {\scriptstyle{{\bf H}}}\mathfrak{I}$, (ii)
$H \in {\scriptstyle{{\bf LH}}}\mathfrak{I}$, (iii)
$H \in \Phi\mathfrak{I}$ and (iv) $H \in \Phi_{flat}\mathfrak{I}$.

(i) Assume that $H \in {\scriptstyle{{\bf H}}}\mathfrak{I}$
and use induction on the ordinal $\alpha$, which is such that
$H \in {\scriptstyle{{\bf H}}}_{\alpha}\mathfrak{I}$. If
$\alpha = 0$, then $H \in \mathfrak{I}$ and hence $N$ is
projective, in view of the very definition of
$\mathfrak{I}$-groups. Let $\alpha >0$ and assume that the
result is known for all
${\scriptstyle{{\bf H}}}_{\beta}\mathfrak{I}$-subgroups
$K \subseteq G$ and all ordinals $\beta$ with $\beta < \alpha$.
The group $H$ acts cellularly on a finite dimensional
contractible CW-complex $X$, with each cell stabilizer
contained in ${\scriptstyle{{\bf H}}}_{\beta}\mathfrak{C}$
for a suitable $\beta < \alpha$. Then, there is an induced
short exact sequence of $kH$-modules
\[ 0 \longrightarrow N_d \longrightarrow \cdots
     \longrightarrow N_1 \longrightarrow N_0
     \longrightarrow N \longrightarrow 0 , \]
where $d$ is the dimension of $X$ and each $N_i$ is a direct
sum of modules of the form $\mbox{ind}_K^H\mbox{res}_K^HN$,
for a suitable
${\scriptstyle{{\bf H}}}_{\beta}\mathfrak{I}$-subgroup
$K \subseteq H$ for some $\beta < \alpha$. For any such $K$
we note that
$\mbox{res}_K^HN = \mbox{res}_K^HM \otimes_k \mbox{res}_K^GB(G)$
and $\mbox{res}_K^HM \in \mathcal{P}_0(kK)$. Invoking our
induction hypothesis, we therefore conclude that
$N_i \in {\tt Proj}(kH)$ for all $i=0,1, \ldots ,d$ and hence
$\mbox{pd}_{kH} N \leq d$. This inequality holds for any
$M \in \mathcal{P}_0(kH)$. Since $M$ is the $d$-th syzygy
in a projective resolution of another $kH$-module
$M' \in \mathcal{P}_0(kH)$, a standard argument shows that
the $kH$-module $N$ is projective.

(ii) Assume that $H \in {\scriptstyle{{\bf LH}}}\mathfrak{I}$
and consider a finitely generated subgroup $K \subseteq H$.
Then, there exists an
${\scriptstyle{{\bf H}}}\mathfrak{I}$-subgroup $K' \subseteq H$
with $K \subseteq K'$. Since
$\mbox{res}_{K'}^HN = \mbox{res}_{K'}^HM \otimes_k
 \mbox{res}_{K'}^GB(G)$
and $\mbox{res}_{K'}^HM \in \mathcal{P}_0(kK')$, we conclude
from case (i) above that
$\mbox{res}_{K'}^HN \in {\tt Proj}(kK')$. Restricting further
to $K$, it follows that $\mbox{res}_K^HN \in {\tt Proj}(kK)$.
In particular, the $kK$-module $\mbox{res}_K^HN$ is flat and
hence the induced $kH$-module $\mbox{ind}_K^H\mbox{res}_K^HN$
is flat as well. Since the $kH$-module $N$ is the filtered
colimit of the system $(\mbox{ind}_K^H\mbox{res}_K^HN)_K$,
where $K$ runs through the family of all finitely generated
subgroups of $H$, we conclude that $N$ is flat. This is the
case for any $M \in \mathcal{P}_0(kH)$. Since $M$ is a cokernel
of an acyclic complex of projective $kH$-modules
with all cokernels in $\mathcal{P}_0(kH)$, it follows that
$N$ is a cokernel of an acyclic complex of projective
$kH$-modules, all of whose cokernels are flat. The result
by Benson, Goodearl and Neeman on the contractibility of
pure acyclic complexes of projective modules (cf.\
\cite[Theorem 2.5]{BG} and \cite[Remark 2.15]{N2}) therefore
implies that $N$ is actually a projective $kH$-module.

(iii) Assume that $H \in \Phi\mathfrak{I}$. For any
$\mathfrak{I}$-subgroup $K \subseteq H$ the restricted
$kK$-module
$\mbox{res}_K^HN = \mbox{res}_K^HM \otimes_k \mbox{res}_K^GB(G)$
is projective, in view of the definition of $\mathfrak{I}$-groups,
since $\mbox{res}_K^HM \in \mathcal{P}_0(kK)$; hence,
$\mbox{pd}_{kH}N < \infty$. Since this is the case for any
$M \in \mathcal{P}_0(kH)$ and the class $\mathcal{P}_0(kH)$ is
closed under direct sums, a standard argument shows that there
is an upper bound on the projective dimension of $N$ for all
$M \in \mathcal{P}_0(kH)$. We may now conclude, as in the
last part of the proof of (i) above, that the $kH$-module
$N$ is actually projective.

(iv) Assume that $H \in \Phi_{flat}\mathfrak{I}$. Then,
arguing as in (iii) above, we conclude that for any
$\mathfrak{I}$-subgroup $K \subseteq H$ the $kK$-module
$\mbox{res}_K^HN$ is projective and hence flat; therefore,
$\mbox{fd}_{kH}N < \infty$. Using the same argument as in
(iii), it follows that $N$ is actually flat. Furthermore,
using the result by Benson, Goodearl and Neeman on the
contractibility of pure acyclic complexes of projective
modules (cf.\ \cite[Theorem 2.5]{BG} and \cite[Remark
2.15]{N2}), as in the last part of the proof of (ii)
above, we conclude that $N$ is projective, as needed.
\end{proof}

\noindent
We conclude this section with some consequences of the results
above.

\begin{Corollary}
The closure $\overline{\mathfrak{F}}$ of the class
$\mathfrak{F}$ of finite groups under the operation
${\scriptstyle{{\bf LH}}} \cup \Phi \cup \Phi_{flat}$
is a subclass of $\mathfrak{I}$.
\end{Corollary}

\begin{proof}
This is an immediate consequence of Lemma 5.3(i) and
Theorem 5.4.
\end{proof}

\begin{Corollary}
Let $H$ be an $\overline{\mathfrak{F}}$-subgroup
of a group $G$.
Then:

(i) the induction functor $\mbox{ind}_H^G$ maps
${\tt Cof}(kH)$ into ${\tt Cof}(kG)$ and

(ii) the restriction functor $\mbox{res}_H^G$ maps
${\tt Cof}(kG)^{\perp}$ into ${\tt Cof}(kH)^{\perp}$.
\end{Corollary}

\begin{proof}
This follows from Proposition 5.2, in view of Corollary
5.6.
\end{proof}

\section{An orthogonality criterion for cofibrant modules}

\noindent
Our goal in this section is to obtain a criterion for a
$kG$-module to be right orthogonal to the class
${\tt Cof}(kG)$ of all cofibrant modules, over certain
groups $G$. As an application, we describe conditions
implying that the tensor product of two $kG$-modules be
right orthogonal to the class of cofibrant modules.

We begin with the following variation on a theme by Kendall
\cite{Ke}.

\begin{Proposition}
Let $\mathfrak{D}$ be a class of groups. We consider a group
$G$ and a $kG$-module $M$, such that
$\mbox{res}_H^GM \in {\tt Cof}(kH)^{\perp}$ for any
$\mathfrak{D}$-subgroup $H \subseteq G$.

(i) We have $\mbox{res}_H^GM \in {\tt Cof}(kH)^{\perp}$
for any ${\scriptstyle{{\bf H}}}\mathfrak{D}$-subgroup
$H \subseteq G$.

(ii) If $\mathfrak{D}$ is subgroup-closed, then
$\mbox{res}_H^GM \in {\tt Cof}(kH)^{\perp}$ for any
${\scriptstyle{{\bf LH}}}\mathfrak{D}$-subgroup
$H \subseteq G$.

(iii) If $\mathfrak{D}$ is a subclass of $\mathfrak{I}$,
then $\mbox{res}_H^GM \in {\tt Cof}(kH)^{\perp}$ for any
$\Phi\mathfrak{D}$-subgroup $H \subseteq G$.

(iv) If $\mathfrak{D}$ is a subclass of $\mathfrak{I}$,
then $\mbox{res}_H^GM \in {\tt Cof}(kH)^{\perp}$ for any
$\Phi_{flat}\mathfrak{D}$-subgroup $H \subseteq G$.
\end{Proposition}

\begin{proof}
(i) We shall proceed by induction on the ordinal number $\alpha$,
for which $H \in {\scriptstyle{{\bf H}}}_{\alpha}\mathfrak{D}$.
The case where $\alpha = 0$ is precisely our hypothesis on $M$.
For the inductive step, assume that $\alpha >0$ and the result
is true for all
${\scriptstyle{{\bf H}}}_{\beta}\mathfrak{D}$-subgroups of $G$
and all $\beta < \alpha$. We consider an
${\scriptstyle{{\bf H}}}_{\alpha}\mathfrak{D}$-subgroup
$H \subseteq G$ and let $N \in {\tt Cof}(kH)$. Then, there
exists an exact sequence of $kH$-modules
\[ 0 \longrightarrow N_d \longrightarrow \cdots
     \longrightarrow N_1 \longrightarrow N_0
     \longrightarrow N \longrightarrow 0 , \]
where $d$ is the dimension of the $H$-CW-complex witnessing
that $H \in {\scriptstyle{{\bf H}}}_{\alpha}\mathfrak{D}$ and
each $N_i$ is a direct sum of $kH$-modules of the form
$\mbox{ind}_F^H\mbox{res}_F^HN$, with $F$ an
${\scriptstyle{{\bf H}}}_{\beta}\mathfrak{D}$-subgroup of $H$
for some $\beta < \alpha$. Since
$\mbox{res}_F^HN \in {\tt Cof}(kF)$, our induction hypothesis
on $\mbox{res}_F^GM$ implies that
\[ \mbox{Ext}^n_{kH} \! \left( \mbox{ind}_F^H\mbox{res}_F^HN,
   \mbox{res}_H^GM \right) \! = \mbox{Ext}^n_{kF} \! \left(
   \mbox{res}_F^HN,\mbox{res}_F^GM \right) \! = 0 \]
for such a subgroup $F$ and all $n>0$. It follows that
$\mbox{Ext}_{kH}^n \! \left( N_i,\mbox{res}_H^GM \right) \! =0$
for all $i$ and all $n>0$. Then, using dimension shifting, we
conclude that
$\mbox{Ext}_{kH}^n \! \left( N,\mbox{res}_H^GM \right) \! =0$
for all $n>d$. This holds for any $N \in {\tt Cof}(kH)$. Since
the cofibrant $kH$-module $N$ is the $d$-th syzygy in a projective
resolution of another cofibrant $kH$-module $N'$ (cf.\ Proposition
3.2(iii)), it follows that
\[ \mbox{Ext}_{kH}^1 \! \left( N,\mbox{res}_H^GM \right) \! =
   \mbox{Ext}_{kH}^{d+1} \! \left( N',\mbox{res}_H^GM \right)
   \! = 0 . \]
This shows that $\mbox{res}_H^GM \in {\tt Cof}(kH)^{\perp}$ and
completes the inductive step of the proof.

(ii) We now assume that $\mathfrak{D}$ is subgroup-closed and
consider an
${\scriptstyle{{\bf LH}}}\mathfrak{D}$-subgroup $H$ of $G$.
We proceed by induction on the cardinality $\kappa$ of $H$.
If $\kappa \leq \aleph_0$, then the
${\scriptstyle{{\bf LH}}}\mathfrak{D}$-group $H$ is actually
contained in ${\scriptstyle{{\bf H}}}\mathfrak{D}$ and we are
done by (i) above. If $\kappa$ is uncountable, then we may
express $H$ as a continuous ascending union of subgroups
$(H_{\lambda})_{\lambda < \kappa}$, each one having cardinality
$< \kappa$. Since the class ${\scriptstyle{{\bf LH}}}\mathfrak{D}$
is also subgroup-closed, $H_{\lambda}$ is an
${\scriptstyle{{\bf LH}}}\mathfrak{D}$-group and our induction
hypothesis implies that
$\mbox{res}^G_{H_{\lambda}}M \in {\tt Cof}(kH_{\lambda})^{\perp}$
for all $\lambda$. We fix $N \in {\tt Cof}(kH)$ and note that
$\mbox{res}^H_{H_{\lambda}}N \in {\tt Cof}(kH_{\lambda})$ for
all $\lambda$. Letting
$N_{\lambda} =
 \mbox{ind}_{H_{\lambda}}^H \mbox{res}^H_{H_{\lambda}}N$,
we compute
\[ \mbox{Ext}^n_{kH} \! \left( N_{\lambda},\mbox{res}_H^GM
   \right) \! = \mbox{Ext}^n_{kH_{\lambda}} \! \left(
   \mbox{res}^H_{H_{\lambda}}N,\mbox{res}_{H_{\lambda}}^GM
   \right) \! = 0 \]
for all $\lambda$ and $n>0$. The continuous ascending family
of subgroups $(H_{\lambda})_{\lambda < \kappa}$ induces a
continuous direct system of $kH$-modules
$(N_{\lambda})_{\lambda < \kappa}$ with surjective structure
maps, whose colimit is $N$. The short exact sequence of
$kH$-modules
\[ 0 \longrightarrow K_{\lambda} \longrightarrow N_0
     \longrightarrow N_{\lambda} \longrightarrow 0 , \]
where $K_{\lambda}$ is the kernel of the structure map
$N_0 \longrightarrow N_{\lambda}$, is the $\lambda$-th term
of a continuous direct system of short exact sequences. The
colimit of the latter direct system of short exact sequences
is the short exact sequence of $kH$-modules
\begin{equation}
 0 \longrightarrow K \longrightarrow N_0
   \longrightarrow N \longrightarrow 0 .
\end{equation}
Here, $K$ is equal to the continuous ascending union of its
submodules $(K_{\lambda})_{\lambda < \kappa}$. Since
$K_{\lambda +1}/K_{\lambda}$ can be identified with
the kernel of the (surjective) structure map
$N_{\lambda} \longrightarrow N_{\lambda +1}$, we conclude that
$\mbox{Ext}^n_{kH}(K_{\lambda +1}/K_{\lambda},\mbox{res}_H^GM)=0$
for all $\lambda$ and $n>0$. Therefore, Eklof's lemma
\cite[Theorem 7.3.4]{EJ} implies that
$\mbox{Ext}^n_{kH} \! \left( K,\mbox{res}_H^GM \right) \! = 0$
for all $n>0$. The short exact sequence (5) above then implies
that
$\mbox{Ext}^n_{kH} \! \left( N,\mbox{res}_H^GM \right) \! = 0$
for all $n>1$. This holds for any $N \in {\tt Cof}(kH)$. Since
the cofibrant $kH$-module $N$ is the first syzygy in a projective
resolution of another cofibrant $kH$-module $N'$ (cf.\ Proposition
3.2(iii)), it follows that
\[ \mbox{Ext}_{kH}^1 \! \left( N,\mbox{res}_H^GM \right) \! =
   \mbox{Ext}_{kH}^2 \! \left( N',\mbox{res}_H^GM \right) \! = 0 . \]
This shows that $\mbox{res}_H^GM \in {\tt Cof}(kH)^{\perp}$ and
completes the proof of (ii).

(iii) Let $H$ be a $\Phi\mathfrak{D}$-subgroup of $G$. Invoking
the completeness of the cofibrant cotorsion pair, we obtain a
short exact sequence of $kH$-modules
\[ 0 \longrightarrow \mbox{res}_H^GM \longrightarrow L
     \longrightarrow N \longrightarrow 0 , \]
where $L \in {\tt Cof}(kH)^{\perp}$ and $N \in {\tt Cof}(kH)$.
Then, for any $\mathfrak{D}$-subgroup $F \subseteq H$ we have
an exact sequence of restricted $kF$-modules
\[ 0 \longrightarrow \mbox{res}_F^GM \longrightarrow
     \mbox{res}_F^HL \longrightarrow \mbox{res}_F^HN
     \longrightarrow 0 . \]
Since $F \in \mathfrak{D} \subseteq \mathfrak{I}$, Proposition
5.2(ii) implies that $\mbox{res}_F^HL \in {\tt Cof}(kF)^{\perp}$.
Of course, the restricted $kF$-module $\mbox{res}_F^HN$ is
cofibrant. Since $\mbox{res}_F^GM \in {\tt Cof}(kF)^{\perp}$,
in view of our hypothesis on $M$, the thickness of the class
${\tt Cof}(kF)^{\perp}$ implies that
$\mbox{res}_F^HN \in {\tt Cof}(kF)^{\perp}$. Invoking Proposition
3.3, we conclude that $\mbox{res}_F^HN$ is a projective
$kF$-module. This is the case for any $\mathfrak{D}$-subgroup
$F$ of $H$ and hence $\mbox{pd}_{kH}N < \infty$. It follows
that $N \in {\tt Cof}(kH)^{\perp}$ (cf.\ Proposition 3.5) and
hence the thickness of ${\tt Cof}(kH)^{\perp}$ implies that
$\mbox{res}_H^GM \in {\tt Cof}(kH)^{\perp}$, as needed.

(iv) We fix a $\Phi_{flat}\mathfrak{D}$-subgroup
$H \subseteq G$ and follow the argument in (iii)
above to conclude that the module $N$ therein has finite
flat dimension. Then, Proposition 3.5 implies that
$N \in {\tt Cof}(kH)^{\perp}$. Invoking the thickness
of ${\tt Cof}(kH)^{\perp}$, we conclude that
$\mbox{res}_H^GM \in {\tt Cof}(kH)^{\perp}$.
\end{proof}

\begin{Corollary}
Let $\mathfrak{D}$ be a subgroup-closed subclass of $\mathfrak{I}$
and denote by $\overline{\mathfrak{D}}$ its closure under the
operation ${\scriptstyle{{\bf LH}}} \cup \Phi \cup \Phi_{flat}$.
We consider a group $G$ and a $kG$-module $M$, such that
$\mbox{res}_H^GM \in {\tt Cof}(kH)^{\perp}$ for any
$\mathfrak{D}$-subgroup $H \subseteq G$. Then, we have
$\mbox{res}_H^GM \in {\tt Cof}(kH)^{\perp}$ for any
$\overline{\mathfrak{D}}$-subgroup $H \subseteq G$.
\end{Corollary}

\begin{proof}
We recall that groups contained in $\overline{\mathfrak{D}}$
admit a hierarchical description. More precisely, we define
the classes $\mathfrak{D}_{\alpha}$ for any ordinal number
$\alpha$, as follows:
\newline
(a) $\mathfrak{D}_0 = \mathfrak{D}$,
\newline
(b) $\mathfrak{D}_{\alpha +1} =
     {\scriptstyle{{\bf LH}}}\mathfrak{D}_{\alpha}
     \cup \Phi\mathfrak{D}_{\alpha} \cup
     \Phi_{flat}\mathfrak{D}_{\alpha}$
    for any ordinal $\alpha$ and
\newline
(c) $\mathfrak{D}_{\alpha} =
     \bigcup_{\beta < \alpha}\mathfrak{D}_{\beta}$
    for any limit ordinal $\alpha$.
\newline
Then, $\overline{\mathfrak{D}}$ is the class consisting of
those groups $H$, for which $H \in \mathfrak{D}_{\alpha}$ for
some ordinal number $\alpha$. We shall prove the claim in the
statement using induction on the ordinal $\alpha$, for which
$G \in \mathfrak{D}_{\alpha}$. If $\alpha =0$, then
$H \in \mathfrak{D}_0 = \mathfrak{D}$ and
$\mbox{res}_H^GM \in {\tt Cof}(kH)^{\perp}$, in view of our
hypothesis. For the inductive step of the proof, we assume
that $H$ is a $\mathfrak{D}_{\alpha +1}$-subgroup of $G$ and
the result is known for all $\mathfrak{D}_{\alpha}$-subgroups
of $G$. Since $\mathfrak{D}$ is a subclass of $\mathfrak{I}$,
Theorem 5.4 implies that $\overline{\mathfrak{D}}$ is a subclass
of $\mathfrak{I}$ as well. In particular, $\mathfrak{D}_{\alpha}$
is a subclass of $\mathfrak{I}$. Moreover, $\mathfrak{D}$ is
subgroup-closed and hence $\mathfrak{D}_{\alpha}$ is also
subgroup-closed. Therefore, we may apply Proposition 6.1
to the class $\mathfrak{D}_{\alpha}$: Since
\[ H \in \mathfrak{D}_{\alpha +1} =
     {\scriptstyle{{\bf LH}}}\mathfrak{D}_{\alpha}
     \cup \Phi\mathfrak{D}_{\alpha} \cup
     \Phi_{flat}\mathfrak{D}_{\alpha} , \]
we conclude that $\mbox{res}_H^GM \in {\tt Cof}(kH)^{\perp}$,
as needed.
\end{proof}

\noindent
We are now ready to state and prove the main result of this section.

\begin{Theorem}
Let $\mathfrak{D}$ be a subgroup-closed subclass of $\mathfrak{I}$
and consider its closure $\overline{\mathfrak{D}}$ under the
operation ${\scriptstyle{{\bf LH}}} \cup \Phi \cup \Phi_{flat}$.
If $G$ is a $\overline{\mathfrak{D}}$-group, the following two
conditions are equivalent for a $kG$-module $M$:

(i) $M \in {\tt Cof}(kG)^{\perp}$,

(ii) $\mbox{res}_H^GM \in {\tt Cof}(kH)^{\perp}$ for any
$\mathfrak{D}$-subgroup $H \subseteq G$.
\end{Theorem}

\begin{proof}
Since $\mathfrak{D} \subseteq \mathfrak{I}$, the implication
(i)$\rightarrow$(ii) follows from Proposition 5.2(ii). For the
reverse implication, assume that (ii) holds. Then, Proposition
6.2 implies that $\mbox{res}_H^GM \in {\tt Cof}(kH)^{\perp}$ for
any $\overline{\mathfrak{D}}$-subgroup $H \subseteq G$. Since
$G \in \overline{\mathfrak{D}}$, we may let $H=G$ and conclude
that $M \in {\tt Cof}(kG)^{\perp}$, as needed.
\end{proof}

\noindent
Since finite groups are contained in $\mathfrak{I}$ (cf.\
Lemma 5.3(i)), we may specialize Theorem 6.3 to the case
where $\mathfrak{D} = \mathfrak{F}$.

\begin{Corollary}
Let $\overline{\mathfrak{F}}$ be the closure of the class
$\mathfrak{F}$ of finite groups under the operation
${\scriptstyle{{\bf LH}}} \cup \Phi \cup \Phi_{flat}$.
If $G$ is an $\overline{\mathfrak{F}}$-group, the following
two conditions are equivalent for a $kG$-module $M$:

(i) $M \in {\tt Cof}(kG)^{\perp}$,

(ii) $\mbox{res}_H^GM \in {\tt Cof}(kH)^{\perp}$ for any
finite subgroup $H \subseteq G$.
\end{Corollary}

\section{A monoidal model structure in $kG$-Mod}

\noindent
We may interpret the existence of the projective complete and
hereditary cotorsion pair defined by the cofibrant modules in
the language of model categories. Our goal in this section is
to examine the behaviour of the associated model structure in
the category of $kG$-modules with respect to the formation of
tensor products.

Invoking Theorem \ref{thm:cp} and Propositions \ref{prop:Proj}
and \ref{prop:cof+}, we conclude that the triple
\[ \left( {\tt Cof}(kG),{\tt Cof}(kG)^{\perp}, kG\mbox{-Mod} \right) \]
is a Hovey triple; it therefore induces a model structure
$\mathcal{M}$ in the category of $kG$-modules. The exact
subcategory ${\tt Cof}(kG)$ of $kG$-Mod is Frobenius with
projective-injective objects the projective modules; this
follows from \cite[Proposition 5.2]{Gil}, but it can be
easily verified directly as well. Then, the associated stable
category $\underline{\tt Cof}(kG)$ is equivalent to the
homotopy category $\mbox{Ho}(\mathcal{M})$; cf.\
\cite[Proposition 4.4 and Corollary 4.8]{Gil}. More explicitly,
the model structure $\mathcal{M}$ is defined as follows:
\begin{itemize}
\item All $kG$-modules are fibrant, whereas the cofibrant objects
      are the cofibrant $kG$-modules and the trivial objects are
      the modules in the orthogonal class ${\tt Cof}(kG)^{\perp}$.
\item The fibrations (resp.\ the trivial fibrations) are the
      surjective $kG$-linear maps (resp.\ the surjective
      $kG$-linear maps with kernel in ${\tt Cof}(kG)^{\perp}$).
\item The cofibrations (resp.\ the trivial cofibrations) are the
      injective $kG$-linear maps whose cokernel is a cofibrant
      module (resp.\ the injective $kG$-linear maps whose
      cokernel is projective).
\item The weak equivalences are those $kG$-linear maps that may
      be expressed as compositions of a trivial cofibration
      followed by a trivial fibration.
\end{itemize}

\vspace{0.1in}

\noindent
{\bf Remark 7.1.}
The model structure $\mathcal{M}$ defined above is different than
that defined by Benson in \cite[Theorem 10.6]{Ben}. In the special
case where $k$ is Noetherian, Benson's model structure is defined
on the category of countably presented $kG$-modules with finite
cofibrant dimension.
\addtocounter{Lemma}{1}

\vspace{0.1in}

\noindent
The tensor product of $kG$-modules (with the diagonal action
of $G$) endows the category $kG$-Mod with a symmetric monoidal
structure. We shall use \cite[Theorem 4.3.2]{Hov1} and examine
whether the derived tensor product endows the homotopy category
$\mbox{Ho}(\mathcal{M})$ with the structure of a symmetric monoidal
category as well. To that end, it suffices to verify that the model
structure $\mathcal{M}$ is {\em monoidal}, in the sense of
\cite[Definition 4.2.6]{Hov1}. Invoking \cite[$\S $7]{Hov2},
we have to verify that the following conditions are satisfied:

(a) All cofibrations are $k$-pure monomorphisms.

(b) The tensor product of cofibrant $kG$-modules is cofibrant.

(c) The class of projective $kG$-modules is invariant under
tensoring with cofibrant modules.

(d) If $M$ is a cofibrant module, then a special
${\tt Cof}(kG)$-precover $f : N \longrightarrow k$
induces a weak equivalence
$1 \otimes f : M \otimes_k N \longrightarrow M \otimes_k k = M$.
\newline
Assertion (d) is precisely condition 2 in \cite[Definition 4.2.6]{Hov1},
whereas (a), (b) and (c) together imply condition 1 therein; for more
details, see the proof of \cite[Theorem 7.2]{Hov2}.

As far as the validity of these conditions is concerned, we
observe that (a) follows from the definition of cofibrations, since
any cofibrant module is $k$-projective (cf.\ Lemma 3.4(i)) and
hence $k$-flat. Condition (b) is precisely Lemma 3.4(iii), whereas
(c) follows since the class of projective $kG$-modules is invariant
under tensoring with any $k$-projective $kG$-module. It only remains
to examine whether condition (d) holds. We shall prove that this is
indeed the case for any $\overline{\mathfrak{F}}$-group $G$.

\begin{Lemma}
Let $G$ be a finite group and consider two $k$-projective
$kG$-modules $M,N$. If one of them is contained in
${\tt Cof}(kG)^{\perp}$, then
$M \otimes_k N \in {\tt Cof}(kG)^{\perp}$.
\end{Lemma}

\begin{proof}
Assume that $N \in {\tt Cof}(kG)^{\perp}$. Since $N$ is
$k$-projective, the $kG$-module $N \otimes_k B = N \otimes_k kG$
is projective and hence $N$ is cofibrant. On the other hand,
the intersection ${\tt Cof}(kG) \cap {\tt Cof}(kG)^{\perp}$
coincides with the class of projective modules (cf.\
Proposition 3.3), so that $N \in {\tt Proj}(kG)$. Then, the
$kG$-module $M \otimes_k N$ is projective as well and hence
$M \otimes_k N \in {\tt Cof}(kG)^{\perp}$, in view of Proposition
3.5.
\end{proof}

\begin{Corollary}
Let $G$ be an $\overline{\mathfrak{F}}$-group and consider
two $k$-projective $kG$-modules $M,N$. If one of them is
contained in ${\tt Cof}(kG)^{\perp}$, then
$M \otimes_k N \in {\tt Cof}(kG)^{\perp}$.
\end{Corollary}

\begin{proof}
Let $H \subseteq G$ be a finite subgroup and consider
the restricted $kH$-modules $M' = \mbox{res}_H^GM$ and
$N' = \mbox{res}_H^GN$. These $kH$-modules are both
$k$-projective and one of them is contained in
${\tt Cof}(kH)^{\perp}$; cf.\ Corollary 6.4. It follows
from Lemma 7.2 that the $kH$-module
$\mbox{res}_H^G(M \otimes_k N) = M' \otimes_k N'$ is also
contained in ${\tt Cof}(kH)^{\perp}$. Since this is the
case for any finite subgroup $H$ of $G$, we may invoke
Corollary 6.4 once more and conclude that
$M \otimes_k N \in {\tt Cof}(kG)^{\perp}$.
\end{proof}

\noindent
We can now state and prove the main result of this section.

\begin{Theorem}
Let $G$ be an $\overline{\mathfrak{F}}$-group. Then, the
category $kG$-Mod endowed with the model structure $\mathcal{M}$
defined above is a monoidal model category and the homotopy
category $\mbox{Ho}(\mathcal{M})$ is a symmetric monoidal
category.
\end{Theorem}

\begin{proof}
As we have already noted, conditions (a), (b) and (c) that are
stated following Remark 7.1 hold for any group $G$. Regarding
condition (d), let $M$ be a cofibrant $kG$-module and consider
a special ${\tt Cof}(kG)$-precover $f : N \longrightarrow k$.
Then, there is a short exact sequence of $kG$-modules
\begin{equation}
 0 \longrightarrow K \longrightarrow N
   \stackrel{f}{\longrightarrow} k \longrightarrow 0 ,
\end{equation}
where $K = \mbox{ker} \, f \in {\tt Cof}(kG)^{\perp}$. Being
cofibrant, the $kG$-module $N$ is $k$-projective (cf.\ Lemma
3.4(i)). Since the short exact sequence is necessarily $k$-split,
the $kG$-module $K$ is $k$-projective as well. We now consider
the induced short exact sequence of $kG$-modules
\[ 0 \longrightarrow M \otimes_k K \longrightarrow M \otimes_k N
     \stackrel{1 \otimes f}{\longrightarrow} M
     \longrightarrow 0 . \]
The cofibrant $kG$-module $M$ being $k$-projective (cf.\ Lemma
3.4(i)), Proposition 7.3 implies that
$M \otimes_k K \in {\tt Cof}(kG)^{\perp}$. Therefore, $1 \otimes f$
is a trivial fibration and hence a weak equivalence.

Having established condition (d), \cite[$\S $7]{Hov2} implies
that the category $kG$-Mod is a monoidal model category with
respect to the model structure $\mathcal{M}$ and the derived tensor
product endows the homotopy category $\mbox{Ho}(\mathcal{M})$ with
a symmetric monoidal structure. To be more precise, for any
$kG$-module $A$ we choose a special ${\tt Cof}(kG)$-precover
$f_A : Q(A) \longrightarrow A$; the cofibrant $kG$-module $Q(A)$
is a {\em cofibrant replacement} of $A$. In particular, a
cofibrant replacement $Q(k)$ of the trivial $kG$-module $k$ is
the $kG$-module $N$ in (6) above. For two $kG$-modules $A$
and $A'$ with cofibrant replacements $Q(A)$ and $Q(A')$
respectively, the derived tensor product $A \otimes^L A'$ is
defined as $Q(A)\otimes_k Q(A')$. The derived tensor product is
symmetric and associative in the homotopy category
$\mbox{Ho}(\mathcal{M})$. The composition
\[ A \otimes^L Q(k) = Q(A) \otimes_k Q(k) = Q(A) \otimes_k N
     \stackrel{1 \otimes f}{\longrightarrow} Q(A) \otimes_k k = Q(A)
     \stackrel{f_A}{\longrightarrow} A \]
is an isomorphism in $\mbox{Ho}(\mathcal{M})$, since both
$1 \otimes f$ and $f_A$ are weak equivalences. It follows that the
homotopy category $\mbox{Ho}(\mathcal{M})$ is a symmetric monoidal
category with unit object $Q(k)$.
\end{proof}

\begin{Corollary}
Let $G$ be an $\overline{\mathfrak{F}}$-group. Then, the
stable category $\underline{\tt Cof}(kG)$ of cofibrant
modules is a symmetric monoidal category.
\end{Corollary}

\medskip

\noindent
{\em Acknowledgments.}
I.\ Emmanouil was supported by the Hellenic Foundation for
Research and Innovation (H.F.R.I.) under the "3rd Call for
H.F.R.I.\ Research Projects to Support Faculty Members and
Researchers", project number 24921. W.\ Ren was supported by
the Natural Science Foundation of Chongqing, China (No.\
CSTB2025NSCQ-GPX1014).

\bigskip

{\footnotesize \noindent Ioannis Emmanouil,
Department of Mathematics, University of Athens, Athens 15784, Greece \\
E-mail: {\tt emmanoui$\symbol{64}$math.uoa.gr}}

\medskip

{\footnotesize \noindent Wei Ren,
 School of Mathematical Sciences, Chongqing Normal University, Chongqing 401331, PR China\\
 E-mail: {\tt wren$\symbol{64}$cqnu.edu.cn}}

\end{document}